\documentclass{article}

\pdfoutput=1

\makeatletter

\@addtoreset{equation}{section}
\makeatother

\input xy
\xyoption{all}

\usepackage[cp1252]{inputenc}
\usepackage[OT1]{fontenc}
\usepackage[english]{babel}

\usepackage{amssymb, amsmath, amsthm,
amscd,
graphicx,
cite,
}

\def\R{{\mathbb R}}
\def\N{{\mathbb N}}
\def\Z{{\mathbb Z}}

\def\CL{\mathcal{L}}
\def\CM{\mathcal{M}}
\def\CU{\mathcal{U}}

\def\ds{\displaystyle}

\def\p{\psi}
\def\a{\alpha}
\def\b{\beta}

\def\D{\Delta}
\def\G{\Gamma}
\def\g{\gamma}
\def\eps{\varepsilon}
\def\f{\varphi}
\def\S{\Sigma}
\def\s{\sigma}
\def\t{\theta}

\def\lam{\lambda}
\def\Lam{\Lambda}

\def\dlam{\dot \lam}
\def\dx{\dot x}
\def\dy{\dot y}

\def\dq{\dot q}
\def\dh{\dot h}

\def\tcut{t_{\operatorname{cut}}}
\def\tt{\mathbf t}

\newcommand{\sn}{\operatorname{sn}\nolimits}
\newcommand{\cn}{\operatorname{cn}\nolimits}
\newcommand{\dn}{\operatorname{dn}\nolimits}
\newcommand{\sgn}{\operatorname{sgn}\nolimits}

\newcommand{\E}{\operatorname{E}\nolimits}

\newcommand{\Exp}{\operatorname{Exp}\nolimits}

\newcommand{\codim}{\operatorname{codim}\nolimits}
\newcommand{\Imm}{\operatorname{Im}\nolimits}
\newcommand{\Hess}{\operatorname{Hess}\nolimits}
\newcommand{\Ker}{\operatorname{Ker}\nolimits}
\newcommand{\Coker}{\operatorname{Coker}\nolimits}
\newcommand{\ind}{\operatorname{ind}\nolimits}
\newcommand{\dgn}{\operatorname{dgn}\nolimits}
\newcommand{\corank}{\operatorname{corank}\nolimits}

\def\skew{\angle}
\def\trans{\pitchfork}
\def\uLam{\underline{\dot{\Lam}}}

\def\vh{\vec h}
\def\vH{\vec H}

\newcommand{\pder}[2]{\frac{\partial \, #1}{\partial \, #2} }
\newcommand{\pdder}[2]{\frac{\partial^2 #1}{\partial \, {#2}^2} }
\newcommand{\be}[1]{\begin{equation}\label{#1}}
\newcommand{\ee}{\end{equation}}
\newcommand{\ddef}[1]{{\em #1\/}}
\newcommand{\map}[3]{#1 \, : \, #2 \to #3}
\newcommand{\mapto}[3]{#1 \, : \, #2 \mapsto #3}

\newcommand{\eq}[1]{$(\protect\ref{#1})$}
\newcommand{\restr}[2]{\left. #1 \right|_{#2}}

\newcommand{\hypoth}[2]{\medskip\noindent{#1}{\em #2}\medskip}

\def\tilt{\widetilde{t}}
\def\tu{\widetilde{u}}

\def\hp{\widehat{p}}

\def\tG{\widetilde{\G}}
\def\tk{\widetilde{k}}
\def\hk{\widehat{k}}

\def\bu{\bar{u}}

\def\bt{\bar{t}}
\def\bk{\bar{k}}

\def\bt{\bar{t}}

\def\ts{\,{\sn \tau}\,}
\def\tc{\,{\cn \tau}\,}

\def\ss{\,{\sn  p}\,}
\def\cc{\,{\cn  p}\,}
\def\dd{\,{\dn  p}\,}
\def\tsp{\,{\sn^2 \tau}\,}

\def\ssp{\,{\sn^2  p}\,}
\def\ccp{\,{\cn^2  p}\,}
\def\ddp{\,{\dn^2  p}\,}

\def\cct{\,{\cn^3  p}\,}
\def\ddt{\,{\dn^3  p}\,}

\def\ssf{\,{\sn^4  p}\,}

\def\ddf{\,{\dn^4  p}\,}

\def\Ho{$(\mathbf{H1})$ } 
\def\Ht{$(\mathbf{H2})$ } 
\def\Hth{$(\mathbf{H3})$ } 
\def\Hf{$(\mathbf{H4})$ } 
\def\pconj{p_1^{\operatorname{conj}}}
\def\tconj{t_1^{\operatorname{conj}}}

\def\then{\quad\Rightarrow\quad}
\def\thens{\ \Rightarrow\ }
\def\iff{\quad\Leftrightarrow\quad}

\newcommand{\twofiglabel}[6]
{
\begin{figure}[htbp]
\includegraphics[width=0.47\textwidth]{#1}
\hfill
\includegraphics[width=0.47\textwidth]{#4}
\\
\parbox[t]{0.45\textwidth}{\caption{#2}\label{#3}}
\hfill
\parbox[t]{0.45\textwidth}{\caption{#5}\label{#6}}
\end{figure}
}

\newtheorem{theorem}{Theorem}[section]
\newtheorem{lemma}{Lemma}[section]
\newtheorem{corollary}{Corollary}[section]
\newtheorem{proposition}{Proposition}[section]

\theoremstyle{remark}
\newtheorem*{remark}{Remark}

\title{Conjugate points in Euler's elastic problem%
\footnote{Work supported by the Russian Foundation for Basic Research, project
No.~05-01-00703-a.}
}

\author{
Yu. L. Sachkov \\
Program Systems Institute \\
Russian Academy of Sciences \\
Pereslavl-Zalessky 152020 Russia \\
E-mail: sachkov@sys.botik.ru}

\begin{document}

\maketitle

\begin{abstract}
For the classical Euler's elastic problem, conjugate points are described. Inflectional elasticae admit the first conjugate point between the first and the third inflection points. All the rest elasticae do not have conjugate points. 

\bigskip
\noindent
\textbf{\small Keywords:} 
Euler elastica, optimal control, conjugate points,  exponential mapping

\bigskip
\noindent
\textbf{\small Mathematics Subject Classification:}
49J15, 93B29, 93C10, 74B20, 74K10, 65D07
\end{abstract}


\newpage
\tableofcontents

\newpage
\section{Introduction}
\label{sec:intro}
This work is devoted to the study of the following problem considered by Leonhard Euler~\cite{euler}. Given an elastic rod in the plane with fixed endpoints and tangents at the endpoints, one should determine possible profiles of the rod under the given boundary conditions. Euler's problem can be stated as the following optimal control problem:
\begin{align}
&\dx = \cos \t, \label{pr1} \\
&\dy = \sin \t, \label{pr2} \\
&\dot \t = u, \label{pr3} \\
&q=(x,y,\t) \in M= \R^2_{x,y} \times S^1_{\t}, \qquad u \in \R,  \label{pr4} \\
&q(0) = q_0 = (x_0, y_0, \t_0), \qquad
q(t_1) = q_1 = (x_1, y_1, \t_1),
\qquad t_1 \text{ fixed}, \label{pr5} \\
&J = \frac 12 \int_0^{t_1} u^2(t) \, dt \to \min, \label{pr6}
\end{align}
where the integral $J$ evaluates the elastic energy of the rod.

This paper is an immediate  continuation of our previous work~\cite{el_max}, which contained the following material: history of the problem, description of attainable set, proof of existence and boundedness of optimal controls, parametrization of extremals by Jacobi's functions, description of discrete symmetries and the corresponding Maxwell points. In this work we widely use the notation, definitions, and results of work~\cite{el_max}.

Leonhard Euler described extremal trajectories of problem~\eq{pr1}--\eq{pr6}, their projections to the plane $(x,y)$ being called Euler elasticae. Although, the question of optimality of elasticae remained open. Our aim is to characterize global and local optimality of Euler elasticae. Short segments of elasticae are optimal. 
The main result of the previous work~\cite{el_max} in this direction was an upper bound on cut points, i.e., points where elasticae lose their global optimality. In this work we describe conjugate points along elasticae; we obtain precise bounds for the first conjugate point, where the elasticae lose their local optimality.

Each inflectional elastica contains an infinite number of conjugate points. The first conjugate point occurs between Maxwell points; visually, the first conjugate point is located between the first and the third inflection point of the elastica.

All the rest elasticae do not contain conjugate points. 

Notice that Max Born proved in his thesis~\cite{born} that if an elastic arc is free of inflection points, then it does not contain conjugate points, so in this part we repeated Max Born's result. Although, our method of proof is more flexible, and we believe that it will be useful for the study of conjugate points in other optimal control problems.

This work has the following structure.  In Sec.~\ref{sec:conj_gen} we review some basic facts of the theory of conjugate points along regular extremals of optimal control problems. These facts are rather well known, but are scattered through the literature. The main facts of this theory necessary for us are the following: (1)~An instant $t> 0$ is  a conjugate point iff the exponential mapping for time $t$ is degenerate; (2)~Morse index of the second variation of the endpoint mapping along an extremal is equal to the number of conjugate points with account of their multiplicity; (3)~Morse index is equal to Maslov index of the curve in a Lagrange Grassmanian obtained by linearization of the flow of the Hamiltonian system of Pontryagin Maximum Principle; (4)~Maslov index is invariant under homotopies of extremals provided that their endpoints are not conjugate.
We apply this theory for description of conjugate points in Euler's problem. In Section~\ref{sec:inflect} we obtain estimates for the first conjugate point on inflectional elasticae. Moreover, we improve our result of work~\cite{el_max} on the upper bound of cut time on inflectional elasticae. In Section~\ref{sec:noninflect} we show that all the rest elasticae do not contain conjugate points.    
In Section~\ref{sec:final} we present some final remarks on results obtained in this paper and discuss their possible consequences for future work. 

In this work we use extensively Jacobi's functions, see~\cite{lawden, whit_watson}. We apply the system ``Mathematica''~\cite{math} to carry out complicated calculations and to produce illustrations.

\section[Conjugate points, Morse index, and Maslov index]
{Conjugate points, Morse index, \\ and Maslov index}

\label{sec:conj_gen}

In this section we recall some basic facts from the theory of conjugate points in optimal control problems, see~\cite{notes, sar, cime,  arnold_givental, arnold_maslov}.

\subsection{Optimal control problem and Hamiltonians}
\label{subsec:problem_Hamiltonians}
We consider an optimal control problem of the form
\begin{align}
&\dq = f(q,u), \qquad q \in M, \quad u \in U \subset \R^m, \label{sys} \\
&q(0) = q_0, \qquad q(t_1) = q_1, \qquad t_1 \text{ fixed}, \label{bound} \\
&J^{t_1}[u] = \int_0^{t_1} \f(q(t),u(t)) \, dt \to \min, \label{J}
\end{align}
where $M$ is a  finite-dimensional analytic manifold,  $f(q,u)$ and $\f(q,u)$ are analytic in $(q,u)$  families of vector fields   and   functions  on $M$ depending on the control parameter $u \in U$, and $U$ an open subset of $\R^m$. Admissible controls are $u(\cdot) \in L_{\infty}[0, t_1]$, and admissible trajectories $q(\cdot)$ are Lipschitzian.
Let 
$$
h_u(\lam) = \langle \lam, f(q,u)\rangle - \f(q,u), 
\qquad \lam \in T^*M, \quad q = \pi(\lam) \in M, \quad u \in U,
$$
be the \ddef{normal Hamiltonian of PMP} for the problem~\eq{sys}--\eq{J}. Fix a triple $(\tu(t), \lam_t, q(t))$ consisting of a normal extremal control $\tu(t)$, the corresponding extremal $\lam_t$, and the extremal trajectory $q(t)$ for the problem~\eq{sys}--\eq{J}.

In the sequel we suppose that the following hypothesis holds:

\hypoth{\Ho}{For all $\lam \in T^*M$ and  $u \in U$, the quadratic form $\ds \pdder{h_u}{u}(\lam)$ is negative definite.}

Notice that condition~\Ho implies the strong Legendre condition along an extremal pair $(\tu(t), \lam_t)$:
$$
\restr{\pdder{h_u}{u}}{u=\tu(t)} (\lam_t) (v,v) < - \a |v|^2, 
\qquad t \in [0, t_1], \ v \in \R^m, \ \a > 0,
$$
i.e., the extremal $\lam_t$ is regular~\cite{notes}. 

Moreover, we assume also that the following condition is satisfied:

\hypoth{\Ht}
{For any $\lam \in T^* M$, the function $u \mapsto h_u(\lam)$, $u \in U$, has a maximum point $\bu(\lam) \in U$:
$$
h_{\bu(\lam)}(\lam) = \max_{u \in U} h_u(\lam), \qquad \lam \in T^*M.
$$}%
In terms of work~\cite{cime}, condition~\Ht means that $T^*M$ is a regular domain of the Hamiltonian $h_u(\lam)$.  Condition \Ho means that the function $u \mapsto h_u(\lam)$ has no maximum points in addition to $\bu(\lam)$. At the maximum point $\restr{\pder{h}{u}}{u = \bu(\lam)}(\lam) = 0$ for all $\lam \in T^*M$. By implicit function theorem, the mapping $\lam \mapsto \bu(\lam)$ is analytic. The \ddef{maximized Hamiltonian} $H(\lam) = h_{\bu(\lam)}(\lam)$, $\lam \in T^*M$, is also analytic. The extremal  $\lam_t$  is a trajectory of the corresponding Hamiltonian vector field: $\dot \lam_t = \vH(\lam_t)$, and the extremal control is $\tu(t) = \bu(\lam_t)$.

\subsection{Second variation and its Morse index}

Consider the \ddef{endpoint mapping} for problem~\eq{sys}--\eq{J}:
\be{Ft1u}
\map{F_t}{\CU = L_{\infty}([0, t], U)}{M}, 
\qquad u(\cdot) \mapsto (q_u(t), J^t[u]),
\ee
where $q_u(\cdot)$ is the trajectory of the control system~\eq{sys} with the initial condition $q_u(0) = q_0$ corresponding to the control $u = u(\cdot)$. Since $\tu \in \CU$ is an extremal control, it follows that the differential (\ddef{first variation}) $\map{D_{\tu}}{T_{\tu}\CU}{T_{q_{\tu(t)}}M}$ is degenerate, i.e., not surjective, for all $t \in (0, t_1]$,  see~\cite{notes}.

Introduce one more important hypothesis:

\hypoth{\Hth}
{The extremal control $\tu(\cdot)$ is a corank one critical point of the endpoint mapping $F_t$, i.e., 
$$
\codim \Imm D_{\tu}F_t = 1, \qquad t \in (0, t_1].
$$
}%
Condition \Hth means that there exists a unique, up to a nonzero factor, extremal $\lam_t$ corresponding to the extremal control $\tu(t)$.

For any extremal control $u \in \CU$ there is a well-defined Hessian (\ddef{second variation})~\cite{notes} of the endpoint mapping --- a quadratic mapping
$$
\map{\Hess_u F_t}{\Ker D_u F_t}{\Coker D_u F_t = T_{q_u(t)} M/\Imm D_u F_t}.
$$ 
Condition \Hth means that $\dim \left( T_{q_{\tu}(t)} M / \Imm D_{\tu} F_t\right) = 1$ for all $t \in (0, t_1]$, thus uniquely, up to a positive factor, there is defined a quadratic form
\be{Qt}
\map{Q_t = \lam_t \Hess_{\tu} F_t}{\Ker D_{\tu} F_t}{\R}, 
\qquad t \in (0, t_1],
\ee
the projection of the second variation to the extremal $\lam_t$.

The \ddef{Morse index} of a quadratic form $Q$ defined in a Banach space $\CL$ is the maximal dimension of the negative space of the form $Q$:
$$
\ind Q = \max \{ \dim L \mid L \subset \CL, \ \restr{Q}{L \setminus \{0\}} < 0 \}.
$$
The kernel of a quadratic form $Q(x)$ is the space
$$
\Ker Q = \{ x \in \CL \mid Q(x,y) = 0 \ \forall \ y \in \CL\},
$$
where $Q(x,y)$ is the symmetric bilinear form corresponding to the quadratic form $Q(x)$. A quadratic form is called degenerate if it has a nonzero kernel. The multiplicity of degeneration of a form $Q$ is equal to dimension of its kernel: $\dgn Q = \dim \Ker Q$.

Now we return to the quadratic form $Q_t$ given by~\eq{Qt} --- the second variation of the endpoint mapping for the extremal pair $\tu(t), \lam_t$ of the optimal control problem~\eq{sys}--\eq{J}. We continue the quadratic form $Q_t$ from the space $L_{\infty}$ to the space $L_2$ by continuity, and denote by $K_t$ the closure of the space $\Ker D_{\tu} F_t$ in $L_2[0,t]$.

\begin{proposition}[Propos. 20.2~\cite{notes}, Th. 1 \cite{sar}]
Under hypotheses \Ho, \Hth, the quadratic form $\restr{Q_t}{K_t}$ is positive for  small $t > 0$. In particular, $\ind \restr{Q_t}{K_t}=0$ for  small $t > 0$.
\end{proposition}

An instant $t_* \in (0, t_1]$ is called a \ddef{conjugate time} (for the initial instant $t = 0$) along the extremal $\lam_t$ if the quadratic form $\restr{Q_{t_*}}{K_{t_*}}$ is degenerate. In this case the point $q_u(t_*) = \pi(\lam_{t_*})$ is called \ddef{conjugate} for the initial point  $q_0$ along the extremal trajectory $q_u(\cdot)$.

\begin{proposition}[Th. 1 \cite{sar}]
\label{propos:sar}
Under hypotheses \Ho, \Hth$:$
\begin{itemize}
\item[$(1)$]
conjugate points along the extremal $\lam_t$ are isolated: $0 < t_*^1 < \dots < t_*^N \leq t_1$,
\item[$(2)$]
Morse index of the second variation is expressed by the formula
$$
\ind \restr{Q_t}{K_t} = \sum \{ \dgn Q_{t_*^i} \mid 0 < t_*^i < t\}.
$$
\end{itemize}
\end{proposition}

\ddef{Local optimality} of extremal trajectories is characterized in terms of conjugate points. Speaking about local optimality of extremal trajectories in calculus of variations and optimal control, one distinguishes strong optimality (in the norm of the space $C([0, t_1], M)$), and weak optimality  (in the norm of the space $C^1([0, t_1], M)$). Under hypotheses \Ho--\Hth, normal extremal trajectories lose their local optimality (both strong and weak) at the first conjugate point~\cite{notes}. Thus in the sequel, when speaking about local optimality, we mean both strong and weak optimality.

\begin{proposition}[Propos. 21.2, Th. 21.3 \cite{notes}]
\label{propos:Jacobi}
Let conditions \Ho--\Hth be satisfied. 

\begin{itemize}
\item[$(1)$]
If the interval $(0, t_1]$ does not contain conjugate points, then the extremal trajectory $q(t)$, $t \in [0, t_1]$, is locally optimal.
\item[$(2)$]
If the interval $(0, t_1)$ contains a  conjugate point, then the extremal trajectory $q(t)$, $t \in [0, t_1]$, is not locally optimal.
\end{itemize}
\end{proposition}

\subsection{Exponential mapping}
We will add to hypotheses \Ho--\Hth one more condition:

\hypoth{\Hf}
{All trajectories of the Hamiltonian vector field $\vH(\lam)$, $\lam \in T^*M$, are continued to the segment $t \in [0, t_1]$. 
}

Consider the \ddef{exponential mapping for time $t$}:
$$
\map{\Exp_t}{N = T_{q_0}^* M}{M},
\qquad \Exp_t(\lam) = \pi \circ e^{t \vH}(\lam) = q(t), \qquad t \in [0, t_1].
$$
One can construct a theory of conjugate points in terms of the family of the subspaces
$$
\Lam(t) = e_*^{-t\vH} T_{\lam_t}(T^*_{q(t)}M)  \subset T_{\lam_0} (N),
$$
via linearization of the flow of the Hamiltonian vector field $\vH$ along the extremal~$\lam_t$.

\subsection{Maslov index of a curve in Lagrange Grassmanian}
First we recall some basic facts of symplectic geometry, see details in works~\cite{arnold_givental, cime}. 
Let $\S, \s$ be a \ddef{symplectic space}, i.e., $\S$ is a $2n$-dimensional linear space, and $\s$ is a nondegenerate skew-symmetric bilinear form on $\S$. The \ddef{skew-orthogonal complement} to a subspace $\G \subset \S$ is the subspace $\G^{\skew} = \{ x \in \S \mid \s(x,\G) = 0\}$. Since $\s$ is nondegenerate, it follows that $\dim \G^{\skew} = 2n - \dim \G$. A subspace $\G \subset \S$ is called \ddef{Lagrangian} if $\G = \G^{\skew}$, in this case $\dim \G = n$. The set of all Lagrangian subspaces in $\S$ is called \ddef{Lagrange Grassmanian} and is denoted as $L(\S)$, it is a smooth manifold of dimension $n(n+1)/2$ in the Grassmanian $G_n(\S)$ of all $n$-dimensional subspaces in $\S$.

Fix an element $\Pi \in L(\S)$. Define an open set $\Pi^{\trans} = \{ \Lam \in L(\S) \mid \Lam \cap \Pi = 0 \}$. The subset $\CM_{\Pi} = L(\S) \setminus \Pi^{\trans} = \{ \Lam \in L(\S) \mid \Lam \cap \Pi \neq 0 \}$ is called the train for $\Pi$. The set $\CM_{\Pi}$ is not a smooth submanifold in $L(\S)$, but it is represented as a union of smooth strata:  $\CM_{\Pi} = \cup _{k \geq 1} \CM_{\Pi}^{(k)}$, where $\CM_{\Pi}^{(k)} = \{ \Lam \in L(\S) \mid \dim(\Lam \cap \Pi) =k\}$ is a smooth submanifold of $L(\S)$ of codimension $k(k+1)/2$. 
 
Consider a smooth curve $\Lam(t) \in L(\S)$, $t \in [t_0, t_1]$, i.e., a family of Lagrangian subspaces in $\S$ smoothly depending on $t$. Suppose that $\Lam(t_0), \Lam(t_1) \in \Pi^{\trans}$.  \ddef{Maslov index} $\mu_{\Pi}(\Lam(\cdot))$ of the curve $\Lam(\cdot)$ is the intersection index of this curve  with the set $\CM_{\Pi}$. 

In greater detail, let the curve $\Lam(\cdot)$ do not intersect with $\CM_{\Pi} \setminus \CM_{\Pi}^{(1)}$, this can always be achieved by a small perturbation of this curve. For the smooth hypersurface $\CM_{\Pi}^{(1)} \subset L(\S)$, one can define its coorientation in an invariant way as follows.
Any tangent vector to $L(\S)$ at a point $\Lam \in L(\S)$ can naturally be identified with a certain quadratic form on $\Lam$. Take a tangent vector $\dot \Lam(t) \in T_{\Lam(t)}L(\S)$ to a smooth curve $\Lam(t) \in L(\S)$. Choose a point $x \in \Lam(t)$ of the $n$-dimensional space $\Lam(t) \subset \S$. Choose any smooth curve $\tau \mapsto x(\tau)$ in $\S$ such that $x(\tau) \in \Lam(\tau)$ for all $\tau$, and $x(\tau) = x$. Then the quadratic form $\uLam(t)(x)$, $x \in \Lam(t)$, is defined by the formula $\uLam(t)(x) = \s(x, \dx(t))$. One can show that $\s(x, \dx(t))$ does not depend upon the choice of the curve $x(\tau)$, i.e., one obtains a well-defined quadratic form $\uLam(t)$ on the space $\Lam(t)$. Moreover, the correspondence $\dot\Lam \mapsto \uLam$, $ \dot \Lam \in T_{\Lam} L(\S)$, defines an isomorphism of the tangent space $T_{\Lam} L(\S)$ and the linear space of quadratic forms on $\Lam$, see~\cite{cime}.  

Maslov index $\mu_{\Pi}(\Lam(\cdot))$ is defined as the number of transitions of the curve $\Lam(\cdot)$ from the negative side of the manifold $\CM_{\Pi}^{(1)}$ (i.e., with $\uLam(t) > 0$) minus number of reverse transitions (with $\uLam(t) < 0$), taking into account multiplicity.

The fundamental property of Maslov index is its homotopy invariance~\cite{arnold_maslov}: for any homotopy $\Lam^s(t)$, $t \in [t_0^s, t_1^s]$, $s \in [0, 1]$, such that $\Lam^s(t_0^s) , \Lam^s(t_1^s) \in \Pi^{\trans}$ for all $s \in [0, 1]$, we have $\mu_{\Pi}(\Lam^0(\cdot))  = \mu_{\Pi}(\Lam^1(\cdot))$. This fact is proved in the same way as homotopy invariance of the usual intersection index of a curve with smooth cooriented surface.

For monotone curves in Lagrange Grassmanian $L(\S)$ there is the following way of evaluation of Maslov index.

\begin{proposition}[Cor. I.1 \cite{cime}]
\label{propos:maslov_monotone}
Let $\uLam(t) \leq 0$, $t \in [t_0, t_1]$, and let $\{ t \in [t_0, t_1] \mid \Lam(t) \cap \Pi \neq 0 \}$ be a finite subset of the open interval $(t_0, t_1)$. Then
\be{maslov=sum_dim}
\mu_{\Pi}(\Lam(\cdot)) = - \sum _{t \in (t_0, t_1)} \dim (\Lam(t) \cap \Pi).
\ee
\end{proposition}  

In fact, in Cor.~I.1~\cite{cime}, there is given a statement for a nondecreasing curve $(\uLam(t) \geq 0$), then in the right-hand side of formula~\eq{maslov=sum_dim} the sign minus is absent. As indicated in the remark after Cor.~I.1~\cite{cime}, the passage from nondecreasing curves to nonincreasing ones is obtained by the inversion of direction of time $t \mapsto t_0 + t_1 - t$.

The theory of Maslov index can be used for computation of Morse index for regular extremals in optimal control problems.

\subsection{Morse index and Maslov index}
\label{subsec:morse-maslov}
let $\lam_t$, $t \in [0, t_1]$, be a normal extremal of the optimal control problem~\eq{sys}--\eq{J}, and let the hypotheses \Ho--\Hf be satisfied. Consider the family of quadratic forms $Q_t$ given by~\eq{Qt}.

The extremal $\lam_t$ determines a smooth curve
$$
\Lam(t) = e_*^{-t\vH} T_{\lam_t}(T^*_{q(t)}M) \in L(\S), 
\qquad t \in [0, t_1],
$$   
in the Lagrange Grassmanian $L(\S)$, where $\S = T_{\lam_0} (T^* M)$. The initial point of this curve is the tangent space to the fiber $\Lam(0) = \Pi = T_{\lam_0}(T^*_{q_0}M)$. 
The strong Legendre condition (see \Ho) implies monotone decreasing of the curve $\Lam(t)$: the quadratic forms $\uLam(t) < 0$, $t \in [0, t_1]$, see Lemma~I.4~\cite{cime}, thus its Maslov index can be computed via Propos.~\ref{propos:maslov_monotone}. 

On the other hand, the following important statement establishes relation between Morse index of the second variation $Q_t$ and Maslov index of the curve~$\Lam(t)$.

\begin{proposition}[Th. I.3, Cor. I.2 \cite{cime}]
\label{propos:morse_maslov}
Let hypotheses \Ho--\Hf be satisfied. Then:

\begin{itemize}
\item[$(1)$]
An instant $t \in (0, t_1]$ is a conjugate time iff $\Lam(t) \cap \Pi \neq 0$.
\item[$(2)$]
If $\Lam(t_1) \cap \Pi = 0$, then there exists $\bt > 0$ such that 
$$
\ind \restr{Q_{t_1}}{K_{t_1}} = - \mu_{\Pi}(\restr{\Lam(\cdot)}{[t_0, t_1]}) \qquad \forall \ t_0 \in (0, \bt).
$$
\item[$(3)$] 
If $\{t \in (0, t_1] \mid \Lam(t) \cap \Pi \neq 0 \}$ is a finite subset of the open interval $(0, t_1)$, then
$$
\ind \restr{Q_{t_1}}{K_{t_1}} =
\sum _{t \in (0, t_1)} \dim (\Lam(t) \cap \Lam(0)).
$$
\end{itemize}
\end{proposition}

Item (1) of Propos.~\ref{propos:morse_maslov} implies obviously the following  statement.

\begin{corollary}
\label{cor:geom_conj}
Let hypotheses \Ho--\Hf hold. An instant $t \in (0, t_1)$ is a conjugate time iff the mapping $\Exp_t$ is degenerate.
\end{corollary}
\begin{proof}
The condition $\Lam(t) \cap \Pi \neq 0$ means that $e_*^{t\vH}(\Pi) \cap T_{\lam_t}(T^*_{q(t)}M) \neq 0$, which is equivalent to degeneracy of the mapping $\Exp_t = \pi \circ e^{t \vH}$.
\end{proof}

Due to Propos.~\ref{propos:morse_maslov}, we obtain a statement on homotopy invariance of Maslov index of the second variation.

\begin{proposition}
\label{propos:homot_ind}
Let $(u^s(t), \lam^s_t)$, $t \in [0, t^s_1]$, $s \in [0, 1]$, be a continuous in parameter $s$ family of normal extremal pairs in the optimal control problem~\eq{sys}--\eq{J} satisfying the conditions \Ho--\Hf. Assume that for any $s \in [0, 1]$ the terminal instant $t = t^s_1$ is not a conjugate time along the extremal $\lam^s_t$. Then
\be{indQst1=}
\ind \restr{Q_{t_1^1}}{K_{t_1^1}} = \ind \restr{Q_{t_1^0}}{K_{t_1^0}}.
\ee
\end{proposition} 
\begin{proof}
It follows from continuity and strict monotonicity of the curves $\Lam^s(t) = e_*^{-t\vH} T_{\lam_t^s}(T_{q^s(t)}^* M)$, $q^s(t) = \pi(\lam^s_t)$ that there exists $\bt >0 $ such that $\bt < t_s$ for all $s \in [0, 1]$ and any instant $t \in (0, \bt)$ is not a conjugate time along the extremal~$\lam^s_t$.

According to item (2) of Propos.~\ref{propos:morse_maslov}, we have
\be{indQst1}
\ind \restr{Q_{t_1^s}}{K_{t_1^s}} = - \mu_{\Pi} (\restr{\Lam^s(\cdot)}{[t_0, t^s_1]}), 
\qquad \forall \ t_0 \in (0, \bt) \quad \forall s \in [0, 1].
\ee  
For all $s \in [0, 1]$ we have $\Lam^s(t_0) \cap \Pi = \Lam^s(t_1^s) \cap \Pi = 0$. Then homotopy invariance of Maslov index implies that the function $s \mapsto \mu_{\Pi}(\restr{\Lam^s (\cdot)}{[t_0, t^s_1]})$ is constant at the segment $s \in [0, 1]$. Thus equality~\eq{indQst1} implies the required equality~\eq{indQst1=}.
\end{proof}

The following statements can be useful for the proof of absence of conjugate points by homotopy or limit passage.

\begin{corollary}
\label{cor:homot_noconj}
Let all hypotheses of Propos.~$\ref{propos:homot_ind}$ be satisfied. If an extremal trajectory $q^0(t) = \pi(\lam_t^0)$, $t \in (0, t_1^0]$, does not contain conjugate points, then the extremal trajectory $q^1(t) = \pi(\lam_t^1)$, $t \in (0, t_1^1]$, does not contain conjugate points as well.
\end{corollary}
\begin{proof}
A regular extremal does not contain conjugate points iff its Maslov index is zero, so the statement follows from Propos.~\ref{propos:homot_ind}.
\end{proof}

\begin{corollary}
\label{cor:noconj_limit}
Let $(u^s(t), \lam^s_t)$, $t \in [0, + \infty)$, $s \in [0, 1]$, be a continuous in parameter $s$ family of normal extremal pairs in the optimal control problem~\eq{sys}--\eq{J} satisfying the hypotheses \Ho--\Hf. Let for any $s \in [0,1] $ and $T > 0$ the extremal $\lam_t^s$ have no conjugate points for $t \in (0, T]$. Then for any $T > 0$ the extremal $\lam_t^1$ has no conjugate points for $t \in (0, T]$ as well.
\end{corollary}
\begin{proof}
Fix any $T > 0$. By Propos.~\ref{propos:sar}, conjugate points along the extremal $\lam_t^1$ are isolated, thus there exists an instant $t_1 > T$ that is not a conjugate time along $\lam_t^1$. Consider the family of extremals $\lam_t^s$, $t \in [0, t_1]$, $s \in [0, 1]$. Corollary~\ref{cor:homot_noconj} implies that the extremal $\lam_t^1$ has no conjugate points for $t \in (0, t_1]$, thus for $t \in (0, T]$ as well.
\end{proof}

\subsection{Preliminary remarks on  Euler's problem}

In this subsection we show that Euler's elastic problem satisfies all hypotheses required for the general theory of conjugate points described in Subsecs.~\ref{subsec:problem_Hamiltonians}--\ref{subsec:morse-maslov}.
 
Recall~\cite{el_max} that Euler's problem is  stated as follows:
\begin{align}
&\dq = X_1(q) + u X_2(q), \qquad q \in M = \R^2 \times S^1, \quad u \in \R, \label{sys1} \\
&q(0) = q_0, \qquad q(t_1) = q_1, \qquad t_1 \text{ fixed}, \label{boundary1}\\
&J = \frac 12 \int_0^{t_1} u^2 dt \to \min,                  \label{J1}
\end{align}
where $X_1 = \cos \t \pder{}{x} + \sin \t \pder{}{y}$,
$X_2 = \pder{}{\t}$, $[X_1, X_2] = X_3 = \sin \t \pder{}{x} - \cos \t \pder{}{y}$.

This problem has the form~\eq{sys}--\eq{J}, and the regularity conditions for $M$, $f$, $\f$ are satisfied.

In terms of the Hamiltonians $h_i(\lam) = \langle \lam, X_i\rangle$, $\lam \in T^*M$, $i = 1, 2, 3$, the normal Hamiltonian of PMP for Euler's problem is $h_u(\lam) = h_1(\lam) + u h_2(\lam) - \frac 12 u^2$. We have $\ds\pdder{h_u}{u} = - 1 < 0$, i.e., hypothesis \Ho holds.

Condition \Ht obviously holds.

Let $u(t)$ be a normal extremal control in Euler's problem. Corank of the control $u(t)$ is equal to dimension of the space of solutions to the linear Hamiltonian system of PMP $\dlam_t = \vh_1(\lam_t) + u(t) \vh_2(\lam_t)$, i.e., to the number of distinct nonzero solutions to the Hamiltonian system corresponding to the maximized Hamiltonian $H = h_1 + \frac 12 h_2^2$:
\be{Ham_max_sys}
\dlam_t = \vh_1(\lam_t) + h_2 \vh_2(\lam_t), \qquad u(t) = h_2(\lam_t).
\ee
We are interested in the number of distinct nonzero solutions to the vertical subsystem of system~\eq{Ham_max_sys}:
\be{vertHam}
\begin{cases}
\dh_1 = - h_2 h_3, \\
\dh_2 = h_3, \\
\dh_3 = h_1 h_3
\end{cases}
\iff
\begin{cases}
\dot \b = c, \\
\dot c = - r \sin \b, \\
\dot r = 0,
\end{cases}
\ee
where $h_1 = - r \cos \b$, $h_2 = c$, $h_3 = - r \sin \b$, see~\cite{el_max}.

To the extremal control $u(t) \equiv 0$, there correspond two distinct nonzero extremals $(h_1, h_2, h_3)(\lam_t) = (\pm r, 0, 0)$, $r \neq 0$, so in this case $\corank u = 2$.

If $u(t) \not\equiv 0$, then $c_t = u(t) \not\equiv 0$. Then the function $c_t$ determines uniquely via system~\eq{vertHam} the functions $r \sin \b_t = - \dot c_t$ and $r \cos \b_t = -\ddot c_t/c_t$. So the curve $(h_1, h_2, h_3)(\lam_t) \not\equiv 0$ is uniquely determined. Consequently, $\corank u = 1$ in the case $u(t) \not\equiv 0$.

Notice that the control $u(t) \equiv 0$ is optimal, thus in the sequel in the study of optimality of extremal controls we can assume that their corank is equal to~1, i.e., hypothesis \Hth is satisfied.

Finally, hypothesis \Hf is also satisfied since the Hamiltonian field $\vH$ is complete (its trajectories are parametrized by Jacobi's functions determined for all $t \in \R$).

Summing up, all hypotheses \Ho--\Hf are satisfied for Euler's elastic problem, so the theory of conjugate points stated in this section is applicable.

\section{Conjugate points on inflectional elasticae}
\label{sec:inflect}
In this section we describe conjugate points on inflectional elasticae in Euler's  problem. We perform explicit computations and estimates on the basis of parametrization of extremal trajectories obtained in~\cite{el_max}. 

We base upon the decomposition of the preimage of the exponential mapping $T_{q_0}^*M = N = \cup_{i=1}^7 N_i$ introduced in~\cite{el_max}.
In this section we consider the case $\lam \in N_1$. In Subsec.~8.2~\cite{el_max} was obtained a parametrization of the exponential mapping in Euler's problem $\mapto{\Exp_t}{(\f, k, r)}{(x_t,y_t,\t_t)}$ in terms of elliptic coordinates in the domain $N_1$. 
By virtue of Corollary~\ref{cor:geom_conj}, an instant $t$ is a conjugate time iff the mapping $\Exp_t$ is degenerate, i.e., iff its Jacobian $J = \ds\frac{\partial(x_t,y_t,\t_t)}{\partial(\f, k, r)}$ vanishes. A direct computation using parametrization of the exponential mapping obtained in Subsec.~8.2~\cite{el_max}, yields the following:
\begin{align}
&J = \frac{\partial(x_t,y_t,\t_t)}{\partial(\f, k, r)}
=
\frac{1}{\sqrt r \cos(\t_t/2)} \frac{\partial(x_t,y_t,\sin(\t_t/2))}{\partial(\f, k, \sqrt r)} =
- \frac{32 k}{(1-k^2) r^{3/2} \D^2} J_1, 
\label{J=J1}\\
&J_1 = a_0 + a_1 z + a_2 z^2, \qquad z = \tsp \in [0, 1], \label{J1=} \\
&a_2 = - k^2 \ss \, x_1, \label{a2k2sx1} \\
&a_2 + a_1 + a_0 = (1-k^2) \ss \, x_1, \label{a2a1a0x1} \\ 
&a_0 = f_1(p,k) \, x_2, \label{a0=fzx2} \\
&x_1 = - \dd(2 \ss \dd \E^3(p) + ((4k^2 -5) \, p \ss \dd  \nonumber \\
&\qquad + \cc (3 - 6 k^2 \ssp))\E^2(p) + ((4k^2 -5) \cc (1-2k^2 \ssp)\,p \nonumber \\
&\qquad + \ss \dd (4p^2 - 1 + k^2(6\ssp - 4 - 4p^2)))\E(p) \nonumber \\
&\qquad + p \, \ss \dd  (1 -(1-k^2)p^2 + k^2(4k^2-5)\ssp) \nonumber \\ 
&\qquad + 2 \cc(k^2 \ssp \ddp + (1-k^2)(1-2k^2 \ssp)p^2))), \label{x1=} \\
&x_2 = \cc(2(1-k^2)\E(p) - \E^2(p)-(1-k^2)p^2) \nonumber \\
&\qquad +\ss \dd (\E(p)-(1-k^2)p), \label{x2=}  \\
&f_1(p,k) = \ss \dd - (2 \E(p) - p) \cc,   \label{f1=} \\
&p = \sqrt r t /2, \qquad \tau = \sqrt r (\f + t/2), \qquad \D = 1 - k^2 \ssp \tsp. \nonumber   
\end{align}
Here $\cn$, $\sn$, $\dn$, $\E$ are Jacobi's functions, see details in~\cite{el_max}.

\subsection{Preliminary lemmas}
\label{subsec:conjN1_lem}
In this subsection we describe roots and signs of the functions $a_0$ and $a_2 + a_1 + a_0$ that essentially evaluate the numerator of Jacobian $J$ at the extreme points $z = 0$ and 1 respectively, see~\eq{J=J1}, \eq{J1=}.

\subsubsection{Roots of the function $a_0$}
Roots of  the function $f_1(p)$ defined in~\eq{f1=} where described in work~\cite{max3}. For completeness, we cite the statements we will need in the sequel.

\begin{proposition}[Lemma 2.1 \cite{max3}]
\label{propos:lem21max3}
The equation $2E(k) - K(k) = 0$, $k \in [0, 1)$,
has a unique root
$k_0 \in (0, 1)$.
Moreover,
\begin{align*}
&k \in [0, k_0) \then 2E - K > 0, \\
&k \in (k_0, 1) \then 2E - K < 0.
\end{align*}
\end{proposition}

Here and below $K(k)$ and $E(k)$ are complete elliptic integrals of the first and second kinds, see~\cite{el_max, lawden, whit_watson}. Numerical computations give the approximate value $k_0 \approx 0.909$.

\begin{proposition}[Propos. 2.1 \cite{max3}]
\label{propos:prop21max3}
For any
$k \in [0, 1)$
the function
$f_1(p,k)$
has a countable number of roots
$p_n^1$,  $n \in \Z$,
localized as follows: $p^1_0 = 0$ and 
$$
p_n^1 \in (-K + 2 K n, \ K + 2 K n), \qquad n \in \Z.
$$
Moreover, for
$n \in \N$
\begin{align*}
&k \in [0, k_0) \then p_n^1 \in (2K n, K + 2Kn), \\
&k = k_0 \then p_n^1 = 2Kn, \\
&k \in (k_0, 1) \then \ p_n^1 \in (-K + 2Kn, 2Kn),
\end{align*}
where
$k_0$ is the unique root of the equation
$2E(k) - K(k) = 0$,
see Propos.~$\ref{propos:lem21max3}$.
\end{proposition}
 
Now we establish the signs of the function $f_1(p)$ between its zeros $p^1_n$.

\begin{lemma}
\label{lem:fz><0}
For any $m = 0, 1, 2, \dots$, we have: 
\begin{align*}
&p \in (p^1_{2m}, p^1_{2m+1}) \then f_1(p) > 0, \\
&p \in (p^1_{2m+1}, p^1_{2m+2}) \then f_1(p) < 0.
\end{align*}
\end{lemma}
\begin{proof}
By virtue of the equality 
$$
\left(\frac{f_1(p)}{\cc}\right)' = \frac{\ssp \ddp}{\ccp},
$$
the function $f_1(p)/\cc$ increases at the segments of the form $[-K + 2 K n, K + 2 Kn]$, $n \in \Z$. So the function $f_1(p)/\cc$, as well as  $f_1(p)$  changes its sign at the points $p^1_n \in   (-K + 2 K n, K + 2 Kn)$. It remains to verify that $f_1(p)$ is positive at the first interval $(p_0^1, p^1_1) = (0, p^1_1)$. 
We have $f_1(p) = p^3/3 + o(p^3) > 0$, $p \to 0$, and the statement follows.
\end{proof}

Now we describe zeros of the function $x_2$ that enters factorization~\eq{a0=fzx2} of the function $a_0$.

\begin{lemma}
\label{lem:x2=0}
The function $x_2(p)$ given by~\eq{x2=} has a countable number of roots $p = p^{x_2}_n \geq 0$. We have $p^{x_2}_0 = 0$ and $p^{x_2}_n \in (2 K n, K + 2 K n)$ for $n \in \N$, moreover,
\be{pnx2inp1n}
k<k_0 \then p^{x_2}_n \in (p^1_n, K + 2 K n).
\ee
Further,
\begin{align}
&p \in (p^{x_2}_{2m}, p^{x_2}_{2m+1}) \then x_2(p) > 0, \label{x2>0} \\
&p \in (p^{x_2}_{2m+1}, p^{x_2}_{2m+2}) \then x_2(p) < 0, \qquad m = 0, 1, 2, \dots. \label{x2<0} 
\end{align}
\end{lemma}
\begin{proof}
First we show that the function
\be{x2sdup}
\frac{x_2(p)}{\ss \dd} \text{ increases when } p \in (2 Kn, 2 K + 2 Kn).
\ee
Direct computation gives 
\begin{align}
&\left(\frac{x_2(p)}{\ss \dd}\right)' = \frac{x_3(p)}{\ssp \ddp}, \label{x2x3(p)}\\
&x_3 = k^2 (\ccp \E(p) + \a)^2 + (1-k^2)(\E(p) + \b)^2, \label{x3=} \\
&\a = (1-k^2) p \ssp - \cc \ss \dd, \qquad \b = - p \ddp. \nonumber
\end{align}
 Since 
$$
\E(p) + \b = \frac 2 3 k^2 p^3 + o(p^3), 
\qquad 
\ccp \E(p) + \a = \frac 2 3 (1-k^2)p^3 + o(p^3),
$$
then $\E(p) + \b \not\equiv 0$, $\ccp \E(p) + \a \not\equiv 0$. 
So the function $x_3(p)$ given by~\eq{x3=} is nonnegative and vanishes only at isolated points.
By virtue of equality~\eq{x2x3(p)}, assertion~\eq{x2sdup} follows.

Further, we have 
\begin{align*}
&\restr{x_2}{p = 2 Kn} = \cc \, x_4(p), \\
&x_4 = -((1-k^2)(\E(p)-p)^2 + k^2 \E^2(p))<0 \text{ for all } p \neq 0.
\end{align*}
 Thus 
\begin{align*} 
&p = 2 K + 4 K n \then \cc < 0, \  x_2 > 0, \\
&p = 4 Kn \then \cc>0, \  x_2 < 0.
\end{align*}
Consequently, $x_2/(\ss \dd) \to \pm \infty$ as $p \to 2 Kn \mp 0$, $n \in \N$. Moreover, it follows from the asymptotics 
\be{x2=4/45}
x_2 = (4/45)\, k^2(1-k^2)p^6 + o(p^6), \qquad p \to 0,
\ee
that $x_2/(\ss \dd) \to + 0$ as $ p \to + 0$. 

Thus  
\begin{align*}
&p \in (0, 2 K) \then \frac{x_2}{\ss \dd} >  0, \\
&p \in (2 Kn, 2 K + 2 K n) \then \frac{x_2}{\ss \dd} \text{ increases from } -\infty \text{ to } + \infty.
\end{align*}
So there exists a unique root of $x_2(p)/(\ss \dd)$, thus of $x_2(p)$ at the interval $(2 Kn, 2 K + 2 Kn)$, we denote it as $p_n^{x_2}$.

Now we localize $p_n^{x_2}$ w.r.t. the point $K + 2 Kn$. We have 
\begin{align*}
&\restr{x_2}{p = K + 2 K n} = \dd \ss (\E(p) - (1-k^2)p), \\
&\E(p)-(1-k^2)p = k^2 \int_0^p \cn^2 t \, dt > 0, \qquad p > 0.
\end{align*} 
Now 
\begin{align*}
&p = K + 4 Kn \then \ss = 1, \ x_2 > 0 \then \frac{x_2}{\ss} > 0, \\
&p = 3 K + 4 K n \then \ss = -1, \ x_2 < 0 \then \frac{x_2}{\ss} > 0.
\end{align*} 
Consequently, $p_n^{x_2} \in (2 K n, K + 2 K n)$ for all $n \in \N$.

Let $k < k_0$, then $p_n^1 \in (2 K n, K + 2 K n)$, we clarify now the mutual disposition of the points $p^1_n$ and $p^{x_2}_n$ in this case. By virtue of~\eq{f1=}, 
$$
f_1(p) = 0 \iff 
\E(p) = (\dd \ss /\cc + p)/2.
$$
Direct computation gives 
$$
\restr{x_2}{\E(p) =  (\dd \ss /\cc + p)/2} = 8 \ccp (\ss \dd - p \cc) \E(p).
$$ 
Since for $p = p^1_n$ we have $\cc \neq 0$, it follows that  for $p = p^1_n$ the functions $x_2$ and $\ss \dd  - p \cc $ have the same sign. Then both for $p=p^1_{2l-1} \in (4 K l - 2 K, 4 K l - K)$ and for $p=p^1_{2l} \in (4 K l, 4 K l + K)$ we obtain $x_2 /\ss < 0$. Consequently, $p^1_n < p^{x_5}_n$ for all $n \in \N$, i.e., inclusion~\eq{pnx2inp1n} is proved. The roots $p^{x_2}_n$ are localized as required.

For $p>0$ the functions $x_2$ and $\ss$ have distinct roots, so it follows from~\eq{x2sdup} that $x_2$ changes its sign at the points $p^{x_2}_n$, $n \in \N$. The distribution of signs~\eq{x2>0}, \eq{x2<0} follows from the fact that at the first interval $(p^{x_2}_0, p^{x_2}_1) = (0, p^{x_2}_1)$ the function $x_2$ is positive, see~\eq{x2=4/45}.    
\end{proof}

For $p>0$, the function $a_0$ vanishes at the points $p=p_n^1$ and $p=p_n^{x_2}$ defined and localized in Propos.~\ref{propos:prop21max3} and Lemma~\ref{lem:x2=0}.
Now decomposition~\eq{a0=fzx2} and Lemmas~\ref{lem:fz><0}, \ref{lem:x2=0} imply the following statement on distribution of signs of the function $a_0$.

\begin{lemma}
\label{lem:a0><0}
Let $k \in (0, 1)$. If $p \in (0, p_1^1)$, then $a_0> 0$. For any $n \in \N$, if $p \in (p_n^1, p_n^{x_2})$, then $a_0 < 0$, and if $p \in (p_n^{x_2}, p_{n+1}^{1})$, then $a_0 > 0$.
\end{lemma}

\subsubsection{Roots of the function $a_0 + a_1 + a_2$}

In order to obtain a similar description for the function $a_0 + a_1 + a_2$, we have to describe roots of the function $x_1$, see decomposition~\eq{a2a1a0x1}.

\begin{lemma}
\label{lem:x1=0}
For $p \geq 0$, the function $x_1(p)$ defined by~\eq{x1=} has a countable number of roots $p_0 = 0$, $p_n^{x_1} \in (p_n^1, p_{n+1}^1)$, $n \in \N$. Moreover,
\begin{align}
&p \in (p^{x_1}_{2m}, p^{x_1}_{2m+1}) \then x_1(p) > 0, \label{x1>0} \\
&p \in (p^{x_1}_{2m+1}, p^{x_1}_{2m+2}) \then x_1(p) < 0. \label{x1<0} 
\end{align}
\end{lemma}
\begin{proof}
Direct computation gives 
\begin{align}
&\left(\frac{x_1(p)}{\dd f_1(p)}\right)' = \frac{x_5(p)}{4 f_1^2(p)}, \label{x1x5} \\
&x_5 = k^2(\cc E_4 \, p + \a)^2 + (1-k^2)(p E_2 +\b)^2 \geq 0, \nonumber \\
&E_2 = 2 \E(p) - p, \qquad E_4 = \cc (2 \E(p) -p) - 2 \ss \dd, \nonumber\\
&\a = (1+\ssp - 2 k^2 \ssp) E_2^2 + 4 \cc \ss \dd (1-2k^2) E_2 \nonumber\\
&\qquad + 4 (2 k^2 - 1) \ssp \ddp, \nonumber\\
&\b = (2 k^2 \ssp - 1) E_2^2 + 8 k^2 \cc \ss \dd E_2 - 8 k^2 \ssp \ddp. \nonumber
\end{align}
Since 
\begin{align*}
&\cc E_4 \, p + \a = \frac{4}{45} (1-k^2) p^6 + o(p^6) \not\equiv 0, \\
&p E_2 + \b = - \frac{4}{45} k^2 p^6 + o(p^6) \not\equiv 0,
\end{align*}
the function $x_5(p)$ is nonnegative and vanishes at isolated points. In view of equality~\eq{x1x5}, the function $x_1(p)/(\dd f_1(p))$ increases at the intervals where $f_1(p) \neq 0$. 

Now we find the sign of $x_1$ at the points $p^1_n$. We have 
\begin{align*}
&\restr{x_1}{\E(p) = (\dd \ss/\cc + p)/2} = \frac{x_6(p)}{4 \cct}, \\
&x_6(p) = x_6^0 + x_6^1 \, p + x_6^2 \, p^2, \\
&x_6^2(p) = - \ccp \dd(1 - k^2 \ssp(2 - \ssp)), \\
&x_6^1(p) = 2 \cc \ss  (1-2 k^6 \sn^6\, p - k^2 \ssp(3+\ssp) + k^4 \ssf (4 + \ssp)), \\
&x_6^0(p) = - \ddt \ssp (1-k^2 \ssp(2-\ssp)),
\end{align*} 
notice that $1 - k^2 \ssp(2-\ssp) = \ddf + k^2 (1-k^2) \ssf > 0$. Consider the discriminant 
$$
x_{6d} = (x_6^1)^2 - 4 x_6^0 \, x_6^2 = - 16 k^2 (1-k^2) \ccp \sn^6 p \dn^8 p
$$ 
of the quadratic polynomial $x_6(p)$. If $k \neq k_0$, then for $p = p_n^1$ we have $\cc \neq 0$, $\ss \neq 0$, thus $x_6^2, x_6^0, x_{6d}<0$ and $x_6 < 0$. If $k = k_0$, then for $p = p_n^1$ we have $\cc \neq 0$, $\ss \neq 0$, thus $x_6^0 = 0$, $x_6^2 < 0$, $x_{6d} = 0$, $x_6^1 = 0$ and $x_6 = x_6^2 p^2 < 0$. 

So for all $k \in (0, 1)$ if $p = p_n^1 > 0$, then $ \sgn x_1 = - \sgn \cc$. If $p = p^1_{2l-1} \in (4 K l - 2 K, 4 K l - K)$, then $ \cc < 0$, thus $x_1 > 0$. Similarly, if $p = p^1_{2l} \in (4 K l, 4 K l + K)$, then $ \cc > 0$, thus $x_1 < 0$. 

Consequently,
\begin{align}
\label{x1>0p<p11}
&p \in (0, p_1^1) \then \frac{x_1(p)}{\dd f_1(p)} \text{ increases from $0$ to $+ \infty$} \thens x_1(p) > 0, \\
&p \in (p^1_n, p^1_{n+1}), \ n \in \N \then \frac{x_1(p)}{\dd f_1(p)} \text{ increases from } -\infty \text{ to } + \infty, \nonumber
\end{align}
thus $x_1$ has a unique root $p_n^{x_1} \in (p_n^1, p_{n+1}^1)$. 

The required signs of the function $x_1(p)$ at the intervals~\eq{x1>0}, \eq{x1<0} follow from the inequality at the first interval~\eq{x1>0p<p11}, and from the fact that $x_1(p)/(\dd f_1(p))$, thus $x_1(p)$ changes its sign at the points $p_n^{x_1}$, $n \in \N$.  
\end{proof}

\begin{remark}

By virtue of decomposition~\eq{a2a1a0x1}, we have the equality
\begin{align}
\{p>0 \mid a_0 + a_1 + a_2 = 0\} 
&= \{p > 0 \mid \ss = 0 \} \cup \{ p > 0 \mid x_1  = 0 \} \nonumber \\
&= \{ 2 K m \mid m \in \N\} \cup \{ p_n^{x_1} \mid n \in \N\}.
\label{2Kmpx1n}
\end{align}
In order to have a complete description of roots of the function $a_0 + a_1 + a_2$, one should describe mutual disposition of the points $2 K m$ and $p_n^{x_1}$. Numerical computations show that some of these points may coincide one with another. For example, numerical computations yield the following relations between the first roots in the families~\eq{2Kmpx1n}: if $k \in (0, \bk)$, then $p_1^{x_1} > 2 K$; if $k = \bk$, then $ p_1^{x_1} = 2 K$; if $k \in (\bk, 1)$, then $p_1^{x_1} < 2 K$ for a number $\bk \approx 0.998$. We do not go into details of this analysis, but in the sequel we allow different possibilities of mutual disposition of the roots $p_n^{x_1}$ and $2 K m$.
\end{remark}

\subsection{Bounds of conjugate time}
In this subsection we estimate the first conjugate time in Euler's problem along inflectional elasticae.

We obtain from equalities~\eq{a2k2sx1}, \eq{a2a1a0x1} that $a_2 = - k^2/(1-k^2)(a_0 + a_1 + a_2)$, thus the Jacobian appearing in~\eq{J=J1}, \eq{J1=} can be represented as
\be{J1a0a1a2}
J_1(p,k,z) = (1-z) a_0 + z (1-k^2 z)/(1-k^2) (a_0 + a_1 + a_2),
\ee
notice that $(1-k^2 z)/(1-k^2) > 0$. In order to describe the first conjugate point along an extremal trajectory $q(t) = \pi \circ e^{t \vH}(\lam)$, $\lam \in N_1$, it suffices to describe the first positive root of the function $J_1$ for fixed $k$, $z$:
$$
\pconj(k,z) = \min \{ p > 0 \mid J_1(p,k,z) = 0 \}.
$$
This minimum exists since by virtue of regularity of normal extremals, small intervals $p \in (0, \eps)$ do not contain conjugate points. Below in the proof of Th.~\ref{th:pconjN1} we  show this independently on the basis of explicit expression for the function~$J_1$.

\begin{theorem}
\label{th:pconjN1}
Let $\lam \in N_1$. For all $k \in (0, 1)$ and $z \in [0, 1]$, the number $\pconj(k,z)$ belongs to the segment bounded by the points $2 K(k)$, $p_1^1(k)$, namely:
\begin{itemize}
\item[$(1)$]
$k \in (0, k_0) \then \pconj \in [2 K, p_1^1]$,  
\item[$(2)$]
$k = k_0 \then \pconj = 2 K = p_1^1$,  
\item[$(3)$]
$k \in (k_0,1) \then \pconj \in [p_1^1, 2 K]$.
\end{itemize} 
Moreover, for any $k \in (0,1)$ there exists $\eps = \eps(k) > 0$ such that:
\begin{itemize}
\item[$(1')$]
If $k \in (0, k_0)$, then  
\begin{align}
&p \in (0, 2 K) \then J_1 > 0, \label{pin02KJ1>0} \\ 
&p \in (p_1^1, p_1^1 + \eps)  \then J_1 < 0, \label{pinp11J1<0}
\end{align}
\item[$(2')$]
If $k  = k_0$, then 
\begin{align}
&p \in (0, 2 K) \then J_1 > 0,  \label{pin02KJ1>02'} \\ 
&p \in (2K, 2K + \eps) \then J_1 < 0, \label{pin02KJ1<02'} 
\end{align}
\item[$(3')$]
If $k \in (k_0, 1)$, then
\be{pin0p11J1>03'}
p \in (0, p_1^1) \then J_1 > 0,
\ee
moreover,
\begin{itemize}
\item[$(3'a)$]
in the case $p_1^{x_1} \in (p_1^1, 2K):$ 
\be{J1<03'a}
p \in (p_1^{x_1}, p_1^{x_1} + \eps) \then J_1 < 0,
\ee
\item[$(3'b)$]
in the case $p_1^{x_1} = 2K:$ 
\be{J1<=03'b}
p = 2 K = p_1^{x_1}  \then J_1 \leq 0,
\ee
\item[$(3'c)$]
in the case $p_1^{x_1} \in (2K, p_1^2):$ 
\be{J1<03'c}
p \in (2K, 2K + \eps) \then J_1 < 0.
\ee
\end{itemize}
\end{itemize}  
\end{theorem}
\begin{proof}
It is easy to see that by virtue of continuity of the function $J_1(p)$, items $(1')$, $(2')$, $(3')$ imply respectively items $(1)$, $(2)$, $(3)$, so we prove statements  $(1')$, $(2')$, $(3')$.

$(1')$ Fix any $k \in (0, k_0)$, then $ 2 K < p_1^1$, see Propos.~\ref{propos:prop21max3}. 

If $p \in (0, 2 K)$, then Lemmas~\ref{lem:fz><0}, \ref{lem:x2=0}, \ref{lem:x1=0} and decompositions~\eq{a0=fzx2}, \eq{a2a1a0x1} imply the following:
\begin{align*}
&f_1 > 0 \text{ and } x_2 > 0 \then a_0 > 0, \\
&\ss > 0 \text{ and } x_1 > 0 \then a_0 + a_1 + a_2 > 0.
\end{align*}
 Then representation~\eq{J1a0a1a2} gives the inequality $J_1(p,z) > 0$ for all $z \in [0, 1]$ and all $p \in (0, 2 K)$. Implication~\eq{pin02KJ1>0} follows. 

Lemmas~\ref{lem:x2=0}, \ref{lem:x1=0} imply respectively that $p_1^{x_2} \in (p_1^1, 3 K)$, $p_1^{x_1} \in (p_1^1, p_2^1)$. Denote $\hp_1 = \min (p_1^{x_2}, p_1^{x_1}) > p_1^1$. 

If $p \in (p_1^1, \hp_1)$, then we obtain from Lemmas~\ref{lem:fz><0}, \ref{lem:x2=0}, \ref{lem:x1=0} and decompositions~\eq{a0=fzx2}, \eq{a2a1a0x1} the following: 
\begin{align*}
&f_1 > 0 \text{ and } x_2 > 0 \then a_0 < 0, \\
&\ss < 0 \text{ and } x_1 > 0 \then a_0 + a_1 + a_2 < 0.
\end{align*} 
Representation~\eq{J1a0a1a2} implies that $J_1(p,z) < 0$ for all $z \in [0, 1]$ and all $p \in (p_1^1, \hp_1)$, i.e., implication~\eq{pinp11J1<0} is proved for $\eps = \hp_1 - p_1^1> 0$. 

$(2')$ Let $k = k_0$. Similarly to item $(1')$, 
\begin{align*}
&p \in (0, 2 K) \then a_0 >0 \text{ and } a_0 + a_1 + a_2 > 0 \then  J_1 > 0,\\
&p \in (2 K, \hp_1) \then a_0 < 0 \text{ and } a_0 + a_1 + a_2 < 0 \then J_1 < 0,
\end{align*} 
where $\hp_1 = \min(p_1^{x_1}, p_1^{x_2}) > 2 K$. Thus implications~\eq{pin02KJ1>02'}, \eq{pin02KJ1<02'} follow for $\eps = \hp_1 - p_1^{x_2} > 0$. 

$(3')$
Let $k \in (k_0, 1)$, then $p_1^1(k) < 2 K(k)$. 

Let $p \in(0, p_1^1)$. Then we have the following: 
\begin{align*}
&f_1>0 \text{ and } x_2>0 \then a_0 > 0,\\
&\ss > 0 \text{ and } x_1 > 0 \then a_0 + a_1 + a_2 = 0.
\end{align*} 
 Thus $J_1 > 0$, and implication~\eq{pin0p11J1>03'} is proved.

$(3'a)$
Consider the case  $p_1^{x_1} \in (p_1^1, 2 K)$. Let $p \in (p_1^{x_1}, 2 K)$, then: since $f_1 < 0$ and $x_2 > 0$, then $a_ 0 < 0$; since $\ss > 0$ and $x_1 < 0$, then $a_0 + a_1 + a_2 < 0$. Thus $J_1 < 0$, and implication~\eq{J1<03'a} follows for $\eps = 2K - p_1^{x_1} > 0$. In this case
\be{p1conjinp11p1x1}
\pconj(z) \in [p_1^1, p_1^{x_1}] \subset [p_1^1, 2 K) \qquad \forall \ z \in [0, 1].
\ee

$(3'b)$
Consider the case $p_1^{x_1} = 2 K$. Let $p = 2K$, then: since $f_1 < 0$ and $x_2 > 0$, then $a_0 < 0$; since $\ss = x_1 = 0$, then $a_0 + a_1 + a_2 = 0$. Consequently, $J_1 \leq 0$, and implication~\eq{J1<=03'b} follows.

$(3'c)$
Finally, consider the case $p_1^{x_1} \in (2K, p_1^2)$. Let $p \in (2K, \min (p_1^{x_1}, p_1^{x_2}))$, then: since $f_1 < 0$ and $x_2 > 0$, then $a_0 < 0$; since $\ss < 0$ and $x_1 > 0$, then $a_0 + a_1 + a_2 < 0$. Thus $J_1 < 0$, and implication~\eq{J1<03'c} is proved for $\eps = \min(p_1^{x_1}, p_1^{x_2}) - 2 K > 0$.
\end{proof}

\begin{remark}
As one can see from inclusion~\eq{p1conjinp11p1x1}, for $p_1^{x_1} \in (p_1^1, 2 K)$ the range of the function $\pconj(z)$ is strictly less than the segment $[p_1^1, 2 K]$. Judging by the plots of the function $\pconj(z)$, $z = \tsp$, this function is smooth and strictly monotone at the segment $\tau \in [0, K]$, see Figs.~\ref{fig:p1conjk<k0}--\ref{fig:p1conjk>kbar}.

\twofiglabel
{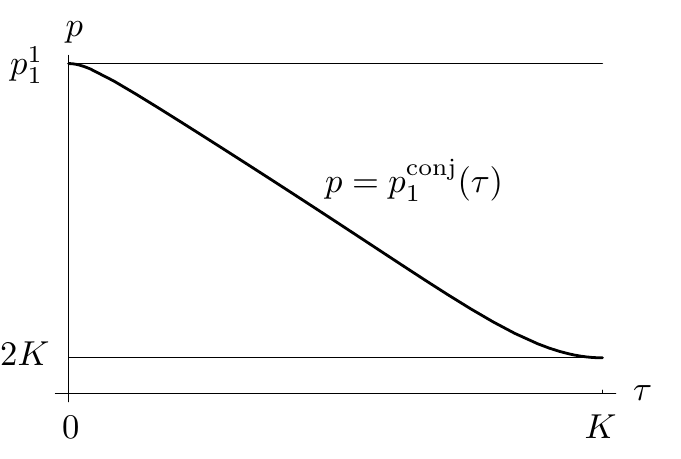}{$p = \pconj(k, \tau)$, $k \in (0, k_0)$}{fig:p1conjk<k0}
{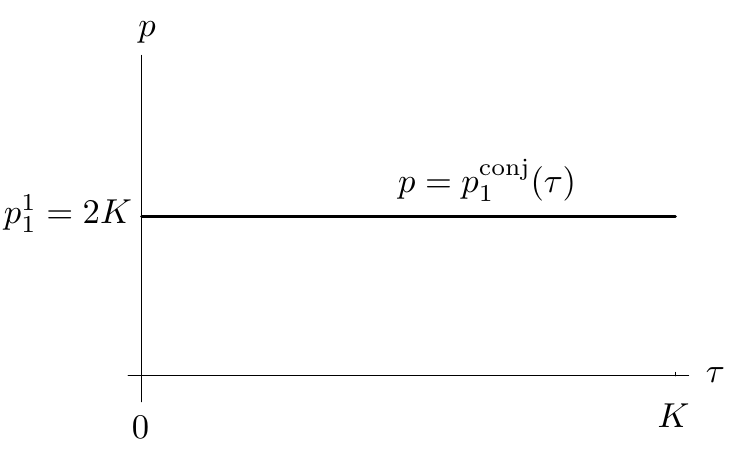}{$p = \pconj(k, \tau)$ $k = k_0$}{fig:p1conjk=k0}


\twofiglabel
{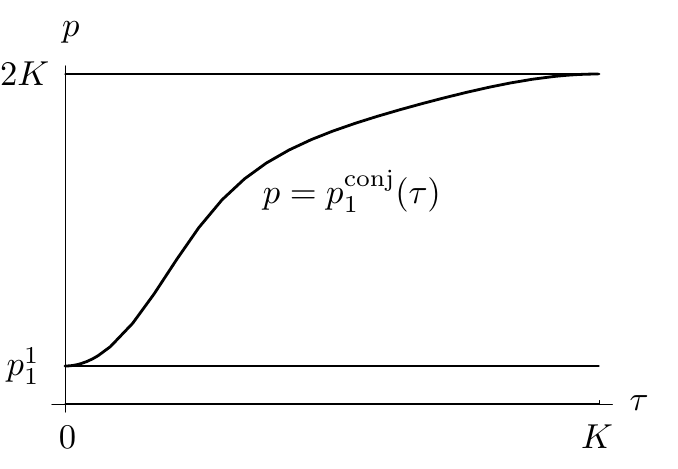}{$p = \pconj(k, \tau)$, $k \in (k_0, 1)$, $2 K \leq p_1^{x_1}$}{fig:p1conjk0<k<kbar}
{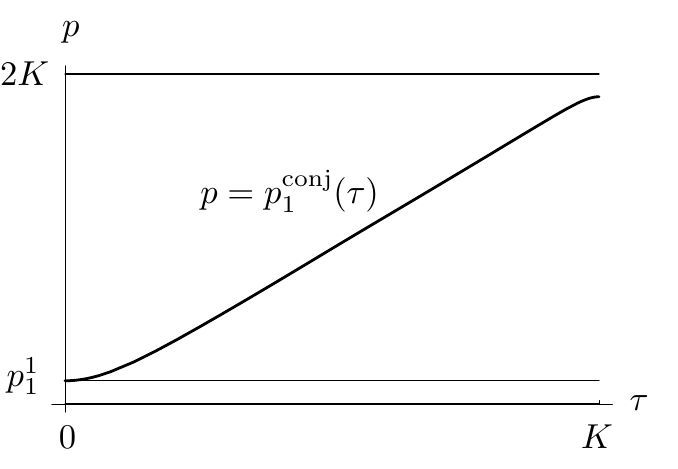}{$p = \pconj(k, \tau)$, $k \in (k_0, 1)$, $2 K > p_1^{x_1}$}{fig:p1conjk>kbar}

\end{remark}

From decompositions~\eq{a0=fzx2}, \eq{a2a1a0x1} and Lemmas~\ref{lem:x2=0}, \ref{lem:x1=0} we obtain the following description of all (not only the first) conjugate points for the cases $z = \tsp = 0$ or 1 (i.e., for elasticae centered respectively at its vertex or inflection point).

\begin{corollary}
\label{cor:pconjz01}
Let $\lam \in N_1$ and $k \in (0,1)$. 
\begin{itemize}
\item[$(1)$]
If $z = 0$, then $\{p > 0 \mid J_1(p,z) = 0\} = \{p_n^1 \mid n \in \N\} \cup \{p_m^{x_2} \mid m \in \N \}$.
\item[$(2)$]
If $z = 1$, then $\{p > 0 \mid J_1(p,z) = 0\} = \{2Kn \mid n \in \N\} \cup \{p_m^{x_1} \mid m \in \N \}$.
\end{itemize}
\end{corollary}

\begin{remark}
According to Lemma~\ref{lem:x2=0} and Propos.~\ref{propos:prop21max3}, in item (1) of Cor.~\ref{cor:pconjz01} all the roots $p_n^1$ and $p_m^{x_2}$ are pairwise distinct. Although, in item (2) some of the roots $2 K n$ and $p_m^{x_1}$ may coincide one with another, see the remark at the end of Subsec.~\ref{subsec:conjN1_lem}. 
\end{remark}

Now we apply preceding results in order to bound the first conjugate time along normal extremal trajectories in the case $\lam \in N_1$:
$$
\tconj(\lam) = \min \{ t > 0 \mid t \text{ conjugate time along trajectory } q(s) = \Exp_s(\lam)\}.
$$

\begin{theorem}
\label{th:tconjN1}
Let $\lam = (k, \f, r) \in N_1$. Then the number $\tconj(\lam)$ belongs to the segment with the endpoints $\ds\frac{4K(k)}{\sqrt r}$, $\ds\frac{2 p_1^1(k)}{\sqrt r}$, namely:
\begin{itemize}
\item[$(1)$]
$k \in (0, k_0) \then \tconj \in \left[\ds\frac{4K(k)}{\sqrt r}, \frac{2 p_1^1(k)}{\sqrt r}\right]$,
\item[$(2)$]
$k = k_0 \then \tconj = \ds\frac{4K(k)}{\sqrt r} = \frac{2 p_1^1(k)}{\sqrt r}$,
\item[$(3)$]
$k \in (k_0,1) \then \tconj \in \left[\ds\frac{2 p_1^1(k)}{\sqrt r}, \frac{4K(k)}{\sqrt r}\right]$.
\end{itemize}
\end{theorem} 
\begin{proof}
By Corollary~\ref{cor:geom_conj},
an instant $t > 0$ is a conjugate time iff
\begin{align*}
&J(t,k,\f,r) = \frac{\partial(x_t, y_t, \t_t)}{\partial(\f, k, r)} = - \frac {32 k}{(1-k^2) r^{3/2} \D^2} J_1(p,k,z) = 0, \\
&p = \sqrt r t /2, \qquad \tau = \sqrt r (\f + t/2), \qquad z = \tsp, \qquad \D = 1 - k^2 \ssp \tsp,
\end{align*}
see~\eq{J=J1}. 

(1) Let $k \in (0, k_0)$, then $\ds\frac{4 K(k)}{\sqrt r} < \frac{2 p_1^1(k)}{\sqrt r}$. According to item $(1')$ of Th.~\ref{th:pconjN1}, for some $\eps = \eps(k) > 0$ we obtain the chains:
$$
t \in \left(0, \ds\frac{4 K}{\sqrt r}\right) \then p \in (0, 2 K) \then J(t,k, \f, r) < 0 \ \forall \ \f, r, 
$$
\begin{align*}
t \in \left(\ds\frac{2 p_1^1(k)}{\sqrt r} , \frac{2 (p_1^1(k)+\eps)}{\sqrt r}\right) &\then p \in (p_1^1(k), p_1^1(k)+ \eps) \\
&\then J(t,k, \f, r) > 0 \ \forall \ \f, r.
\end{align*}
By virtue of continuity of the function $J$ w.r.t. $t$, we obtain the required inclusion $\ds\tconj \in \left[ \frac{4 K(k)}{\sqrt r}, \frac{2 p_1^2(k)}{\sqrt r}\right]$.

Statements (2), (3) of this theorem follow similarly from items $(2')$, $(3')$ of Th.~\ref{th:pconjN1}.
\end{proof}

In Section~12 of work~\cite{el_max}, was defined a function $\map{\tt}{N}{(0, + \infty]}$ that provides an upper bound on cut time in Euler's elastic problem, see Th.~12.1~\cite{el_max}. It follows from formula~(12.2)~\cite{el_max} that
$$
\tt(\lam) = \min\left(\frac{4K(k)}{\sqrt r}, \frac{2 p_1^1(k)}{\sqrt r}\right), \qquad \lam \in N_1.
$$
Comparing this equality with Th.~\ref{th:tconjN1}, we obtain the following statement.

\begin{corollary}
\label{cor:tconjtt}
If $\lam \in N_1$, then $\tconj(\lam) \geq \tt(\lam)$.
\end{corollary}

A natural measure of time along extremal trajectories in Euler's problem is the period of the pendulum $T(k) = 4 K(k)/\sqrt r$. In terms of this measure, the bounds of Th.~\ref{th:tconjN1} are rewritten as follows.

\begin{corollary}
\label{cor:tconjT}
Let $\lam \in N_1$. Then:
\begin{itemize}
\item[$(1)$]
$k \in (0, k_0) \then \tconj \in [T, t_1^1] \subset [T, 3 T/2 ), \quad t_1^1 = 2 p_1^1/\sqrt r \in (T, 3 T/2)$,
\item[$(2)$]
$k = k_0 \then \tconj = T$, 
\item[$(3)$]
$k \in (k_0, 1) \then \tconj \in [t_1^1, T] \subset (T/2 , T], \quad t_1^1 = 2 p_1^1/\sqrt r \in (T/2, T)$.
\end{itemize}
\end{corollary}

It is instructive to state the conditions of local optimality for elastica in terms of their inflection points.

\begin{corollary}
\label{cor:tconj_simple}
Let $\lam \in N_1$, and let $\G = \{\g_s = (x_s, y_s) \mid s \in [0, t]\}$, $q(s) = (x_s, y_s, \t_s) = \Exp(\lam_s)$, be the corresponding elastica.
\begin{itemize}
\item[$(1)$]
If the arc $\G$ does not contain inflection points, then it is locally optimal.
\item[$(2)$]
If  $k \in (0, k_0]$ and the arc $\G$ contains exactly one inflection point, then it is locally optimal.
\item[$(3)$]
If the arc $\G$ contains not less than three inflection points in its interior, then it is not locally optimal.
\end{itemize}
\end{corollary}
\begin{proof}
(1) If the elastic arc $\G$ does not contain inflection points, then its curvature $c_s = 2 k \sqrt r \cn(\sqrt r (\f + s))$ does not vanish for $s \in [0, t]$. But Jacobi's function $\cn(\sqrt r (\f + s))$ vanishes at any segment of length not less than half of its period, thus $t < T/2$. By Cor.~\ref{cor:tconjT}, we have $T/2 < \tconj$, consequently, $t < \tconj$. So the interval $(0, t]$ does not contain conjugate points, thus the corresponding extremal trajectory $q(s)$ is locally optimal (see Propos.~\ref{propos:Jacobi}).

(2) Let $k \in (0, k_0]$, and let the arc $\G$ contain exactly one inflection point. Then the function $c_s$ has exactly one root at the segment $s \in [0, t]$, thus $t < T$. By Cor.~\ref{cor:tconjT}, we have $T \leq \tconj$, so $t < \tconj$, and the elastica $\G$ is locally optimal.

(3)
Let the arc $\G$ contain in its interior not less than 3 inflection points. Then its curvature $c_s$ has not less than 3 roots at the interval $s \in (0, t)$. Consequently, the interval $(0, t)$ contains a complete period $[\tilt_0, \tilt_1]$ of the curvature $c_s$ such that $c_s = 0$ at the endpoints $s = \tilt_0$ and  $s = \tilt_1$, thus $(0,t)$ contains a bigger segment with the same center:
\begin{align*}
&\exists \ [\tilt_0 - \eps, \tilt_1 + \eps] \subset (0, t), \qquad \eps > 0, \\ 
&\sqrt r (\f + \tilt_0) = K + 2 K n, \quad 
\sqrt r (\f + \tilt_1) = 5 K + 2 K n, \qquad n \in \Z.
\end{align*}
So the arc $\G$ contains inside itself the elastica $\tG = \{ \g_s \mid s \in [\tilt_0 - \eps, \tilt_1 + \eps]\}$. Now we show that the arc $\tG$ is not locally optimal, this would mean that the arc $\G$ containing $\tG$ is not locally optimal as well (indeed, if a trajectory $q(s)$, $s \in [0, t]$, is locally optimal, then any its part $q(s)$, $s \in [t_0^1, t_1^1] \subset [0, t]$ is locally optimal as well).

For the arc $\tG$ we have the following: 
\begin{align*}
&(\tilt_1 + \eps) - (\tilt_0 - \eps) = 4 K /\sqrt r + 2 \eps = T + 2 \eps, \\
&\tau = ((\sqrt r (\f + \tilt_0 - \eps) + \sqrt r (\f + \tilt_1 + \eps))/2 = 3 K + 2 Kn, \\
&z = \tsp = 1, \qquad J_1 = a_0 + a_1 + a_2,
\end{align*}
see~\eq{J1a0a1a2}. By Corollary~\ref{cor:pconjz01}, we have $\pconj = \min(2K, p_1^{x_1}) \leq 2 K$, thus $\tconj \leq 4 K/\sqrt r = T$. Consequently, $(\tilt_1 + \eps) - (\tilt_0 - \eps) = T + 2 \eps > \tconj$, and the interval $(\tilt_0 - \eps, \tilt_1 + \eps)$ contains a point $\tconj$ conjugate to the instant $\tilt_0 - \eps$. Thus the arc $\tG$ is not locally optimal, the more so the arc $\G$ is not locally optimal.
\end{proof}

The mathematical notion of local optimality of an extremal trajectory $q(s) = (x_s, y_s, \t_s)$ w.r.t. the functional of elastic energy  corresponds to \ddef{stability} of the corresponding elastica $(x_s,y_s)$.  
Item (3) of Corollary~\ref{cor:tconj_simple} has a simple visual meaning: one cannot keep in hands an elastica having 3 inflection points inside since such an elastica is unstable.

\begin{remark}
In the cases not considered in items (1)--(3) of Corollary~\ref{cor:tconj_simple}, one can find both examples of locally optimal and non-optimal elasticae. 

Let $k > k_0$. If $z = \tsp = 1$ (i.e., the elastica is centered at its inflection point), then by Corollary~\ref{cor:pconjz01}, we have 
$$
\pconj = \min(2K, p_1^{x_1}), \qquad  
p_1^{x_1} \in (p_1^1, p_1^2) \subset (K, 4 K).
$$ 
For $p < K$ we get $p < \pconj$, the corresponding elastica contains one inflection point and is locally optimal, see Fig.~\ref{fig:infl1stable}.
For $p_1^{x_1} <  2 K$ (i.e., for $k \in (\bk, 1)$, $\bk \approx 0.998$) and $p \in (p_1^{x_1}, 2 K)$ we get $p > \pconj = p_1^{x_1}$, the corresponding elastica contains one inflection point and is not locally optimal, see Fig.~\ref{fig:infl1unstable}.

\twofiglabel{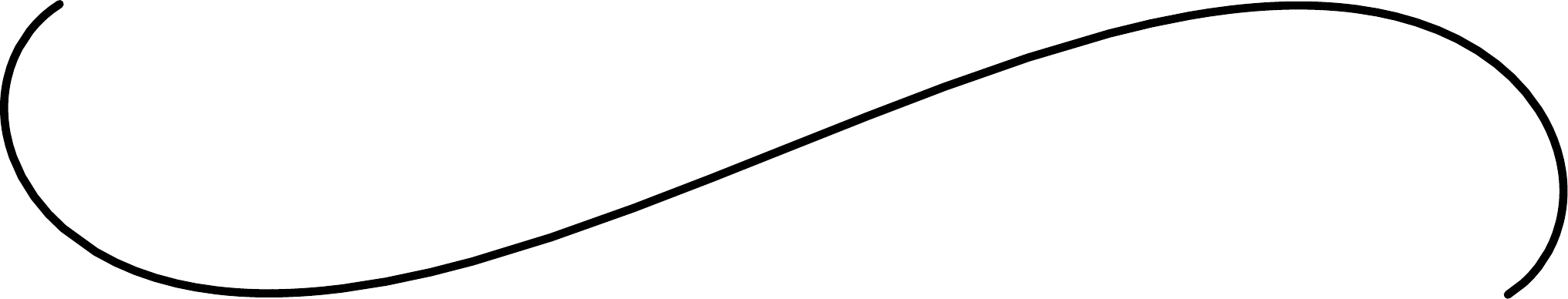}{Locally optimal elastica with 1 inflection point}{fig:infl1stable} 
{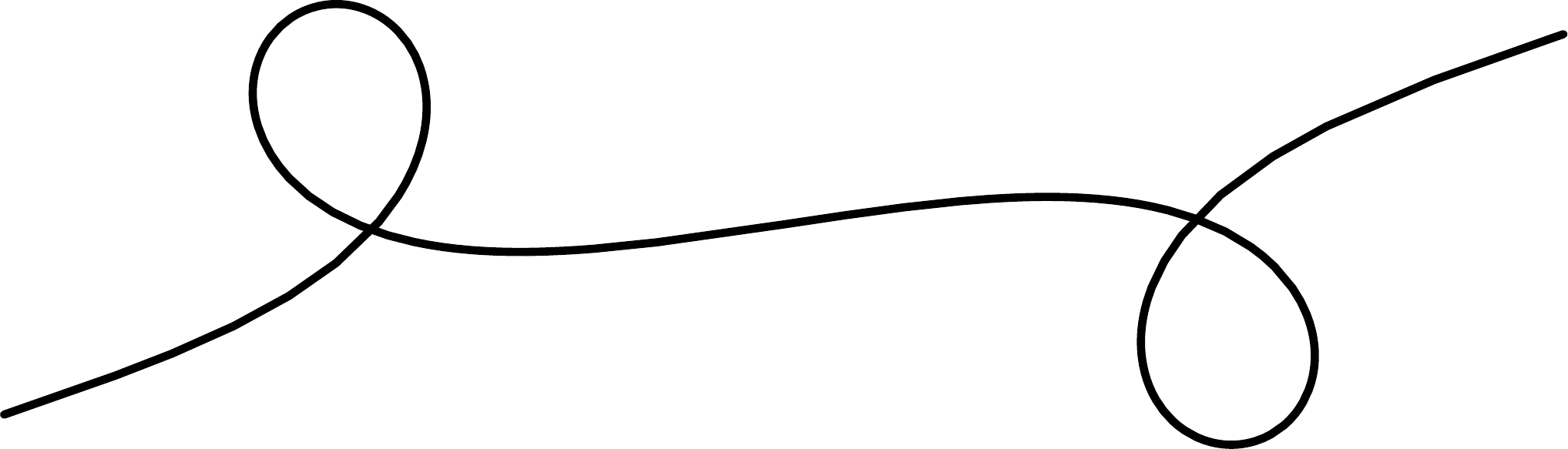}{Locally non-optimal elastica with 1 inflection point}{fig:infl1unstable}

Let $k < k_0$ and $z = \tsp = 0$ (the elastica is centered at its vertex). Then $\pconj = p_1^1 \in (2K, 3 K)$. If $p \in (K, 2 K)$, then $p < p_1^1$, then the corresponding elastica is locally optimal and contains 2 inflection points, see Fig.~\ref{fig:infl2stable}. 

Let $k > k_0$ and $z = \tsp = 0$, then $\pconj = p_1^1 \in (K, 2 K)$. If $p > p_1^1)$, then $p > \pconj$, then the corresponding elastica is not locally optimal and contains 2 inflection points, see Fig.~\ref{fig:infl2unstable}.

\twofiglabel{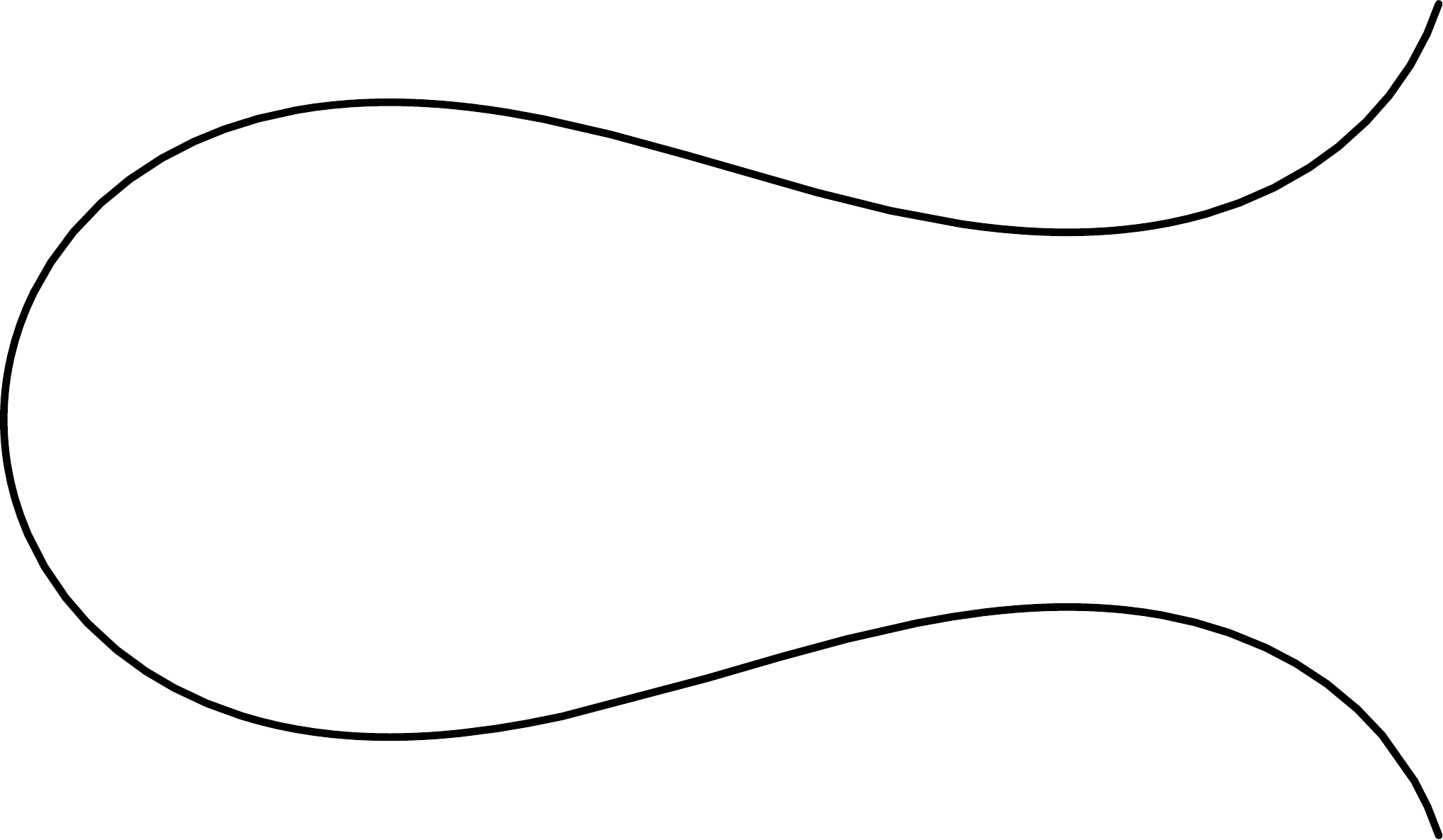}{Locally optimal elastica with 2 inflection points}{fig:infl2stable}
{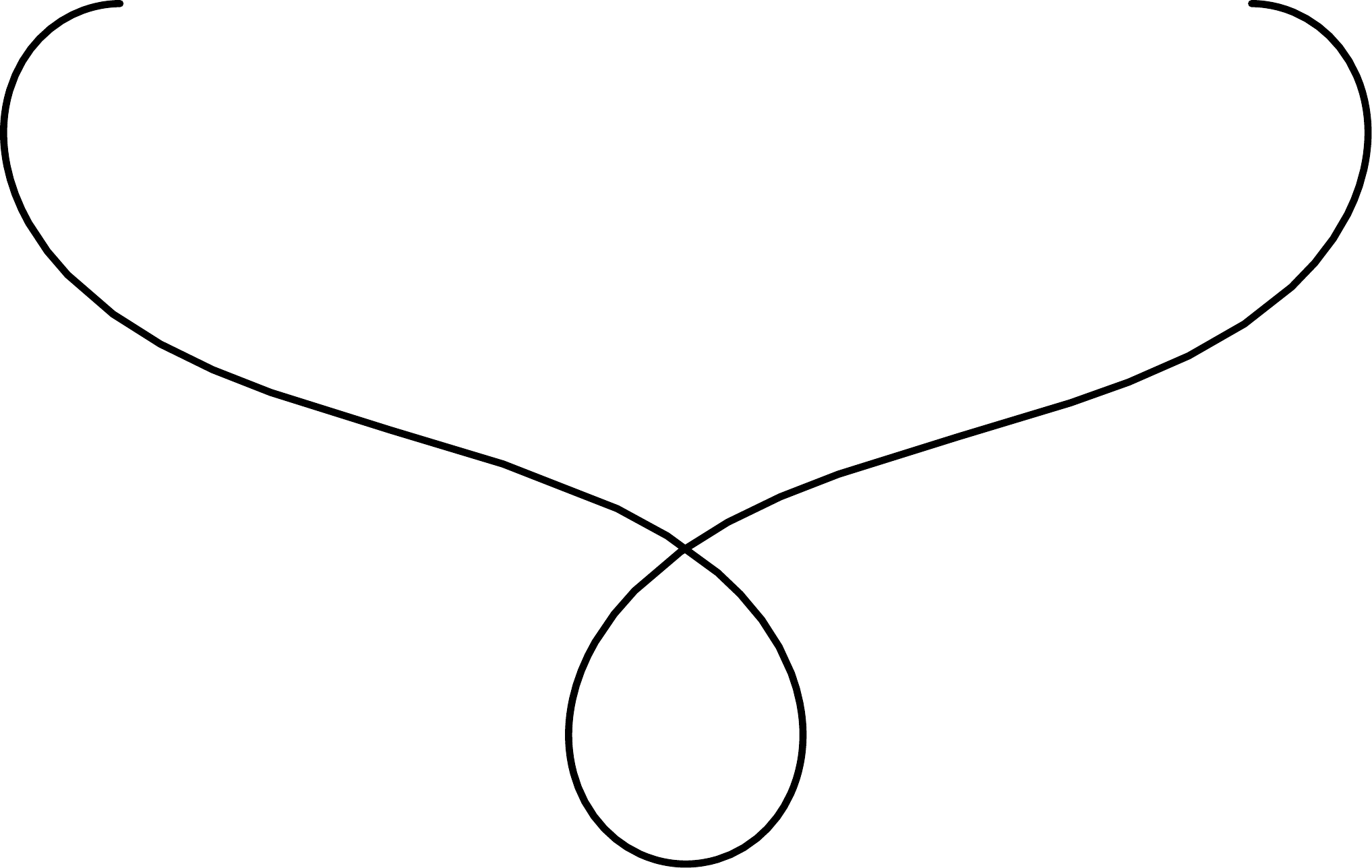}{Locally non-optimal elastica with 2 inflection points}{fig:infl2unstable} 
 
\end{remark}

Corollary~\ref{cor:pconjz01} provides the following description of elastica centered at inflection points or vertices and terminating at conjugate points. 

\begin{corollary}
\label{cor:tconj_tau0K}
Let $\lam  \in N_1$, and let $q(s) = \Exp_s(\lam)$, $s \in [0, t]$, be the corresponding inflectional elastica. 
\begin{itemize}
\item[$(1)$]
If the elastica $q(s)$ is centered at its vertex  (i.e., $\ts = 0$), then the terminal instant $t$ is a conjugate time iff
$$
p = \frac{\sqrt r t}{2 } \in \{p_n^1 \mid n \in \N\} \cup \{p_m^{x_2} \mid m \in \N \}.
$$
\item[$(2)$]
If the elastica $q(s)$ is centered at its inflection point (i.e., $\tc = 0$), then the terminal instant $t$ is a conjugate time iff 
$$
p = \frac{\sqrt r t}{2} \in  \{2Kn \mid n \in \N\} \cup \{p_m^{x_1} \mid m \in \N \}.
$$
 \end{itemize}
\end{corollary}

Figures~\ref{fig:tau0pconj}, \ref{fig:tauKpconj} illustrate respectively cases~(1), (2) of Cor.~\ref{cor:tconj_tau0K}. 

\twofiglabel
{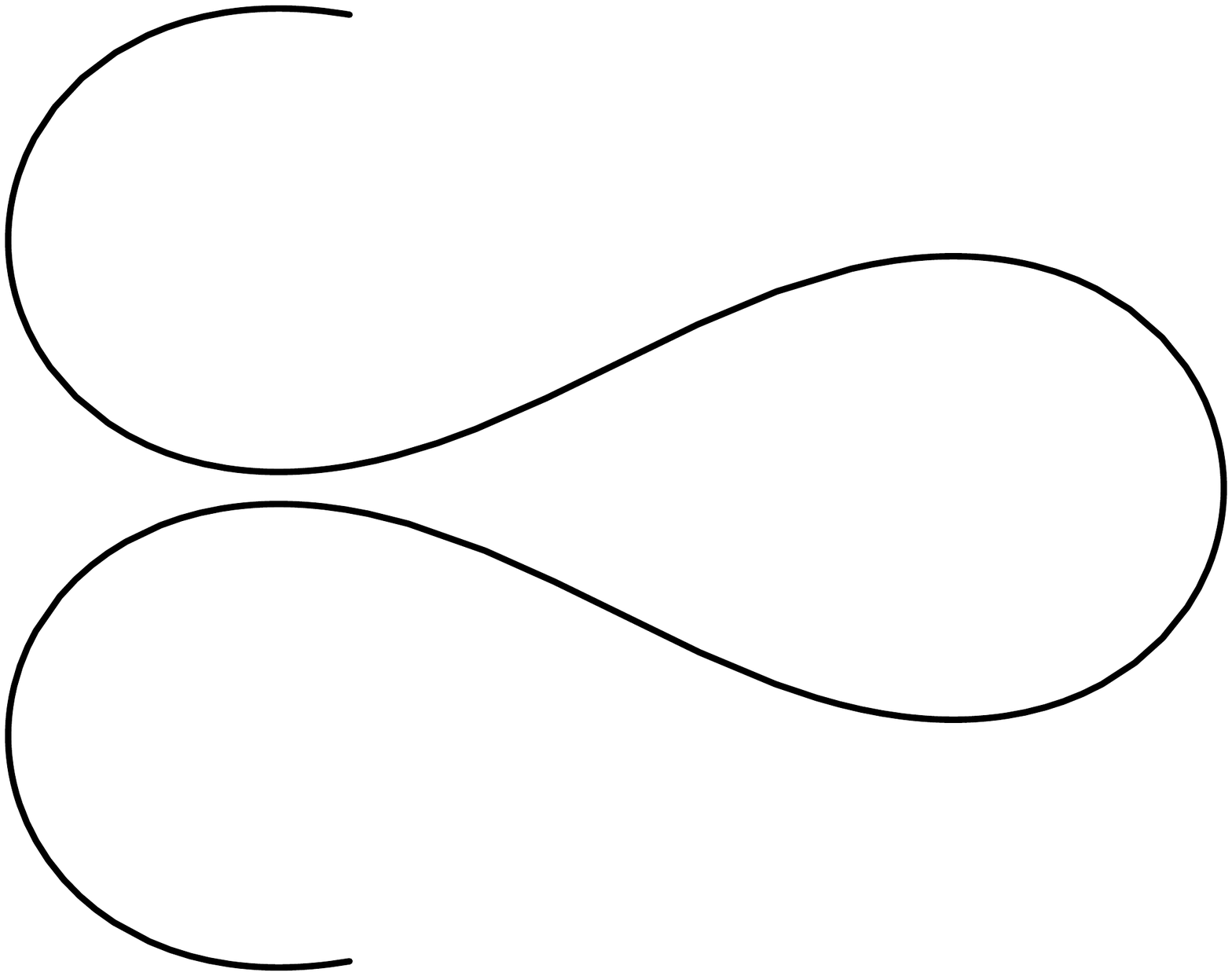}{Conjugate point, $\ts = 0$,  $p = \pconj(k, \tau) = p_1^1(k)$}{fig:tau0pconj}
{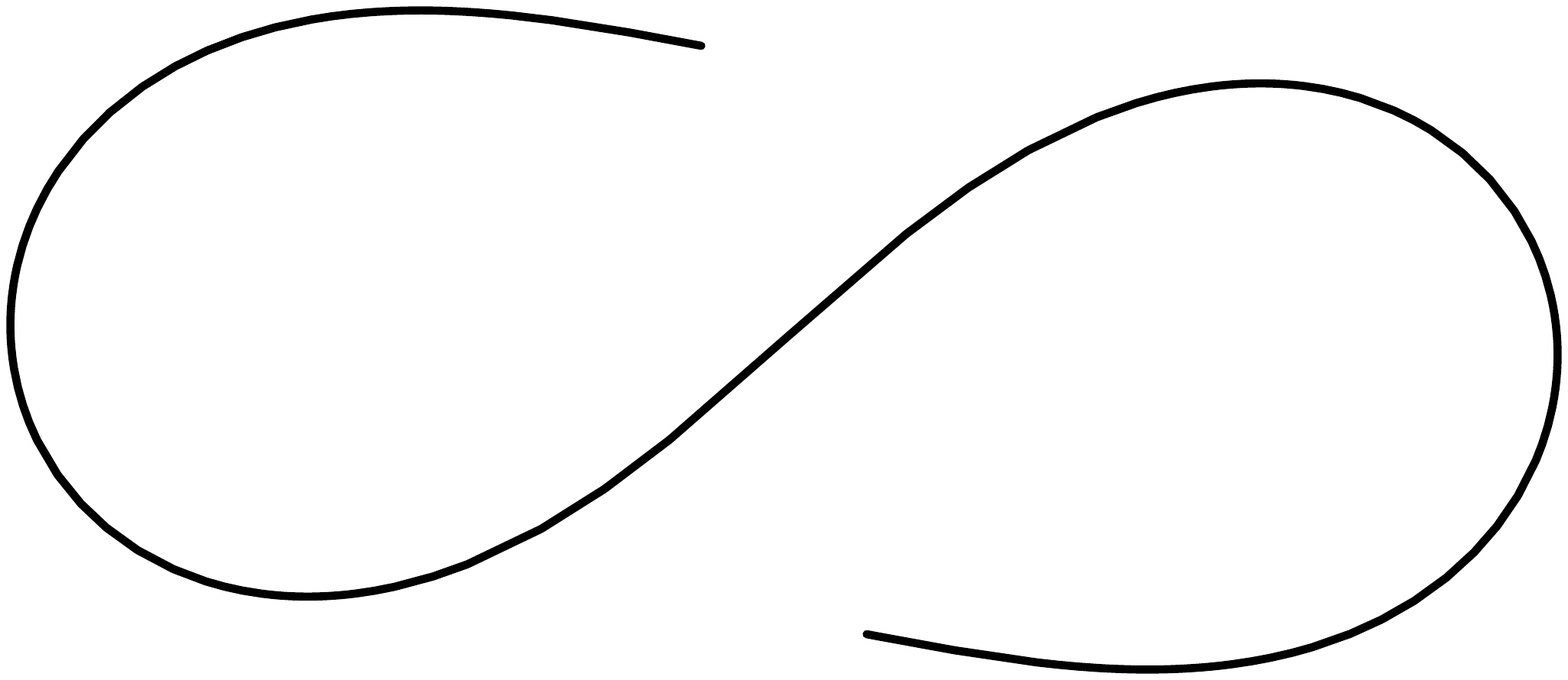}{Conjugate point, $\tc = 0$, $p = \pconj(k, \tau) =  2 K$}{fig:tauKpconj}

\subsection{Upper bound on cut time}

On the basis of results on local optimality obtained in this section, we can improve the statement on  upper bound on the time where elasticae lose their global optimality (i.e., on the  cut time $\tcut(\lam)$), see Th.~12.1~\cite{el_max}. The argument uses the obvious inequality 
$$
\tcut(\lam) \leq \tconj(\lam),
$$   
which follows since if a trajectory is not locally optimal, the more so it is not globally optimal. 

\begin{theorem}
\label{th:tcut_bound}
Let $\lam \in N_1$. Then $\tcut(\lam) \leq \tt(\lam)$.
\end{theorem}
\begin{proof}
We have to prove that the extremal trajectory $q(s) = \Exp_s(\lam)$ is not optimal on any segment of the form $s \in [0, \tt(\lam) + \eps]$, $\eps > 0$. Compute the number $\tau = \frac{\sqrt r}{2}(2 \f + \tt(\lam))$ for the covector $\lam = (k, \f, r)$.

Consider first the case $k \in (0, k_0]$, then $\tt(\lam) = \frac{4K}{\sqrt r}$.
If $\tc \ts \neq 0$, then the equality $\tcut(\lam) \leq \tt(\lam)$ was proved in item (1) of Th.~12.1~\cite{el_max}.
If $\tc = 0$, then the instant $\tt(\lam)$ is a conjugate time by Cor.~\ref{cor:tconj_tau0K}, so the trajectory $q(s)$ is not locally optimal after this instant.
Finally, if $\ts = 0$, then the instant $\tt(\lam)$ is a Maxwell time by item~(1.1) of Th.~11.1~\cite{el_max}.

In the case $k \in (k_0,1)$ we have $\tt(\lam) = \frac{2p_1^1}{\sqrt r}$, and the argument is similar.
If  $\tc \ts \neq 0$, then the statement was proved in  item (1) of Th.~12.1~\cite{el_max}.
If $\ts = 0$, then the instant $\tt(\lam)$ is a conjugate time by  Cor.~\ref{cor:tconj_tau0K}. And if $\tc = 0$, then  the instant $\tt(\lam)$ is a Maxwell time by  item~(1.2) of Th.~11.1~\cite{el_max}.
\end{proof}

\section{Conjugate points on non-inflectional elasticae}
\label{sec:noninflect}

In this section we prove that inflectional elasticae ($\lam \in N_2$), critical elasticae ($\lam \in N_3$), and circles ($\lam \in N_6$) do not contain conjugate points.

Let $\lam \in N_2^+$.  
Similarly to Sec.~\ref{sec:inflect}, we first compute explicitly the Jacobian of the exponential mapping using parametrization of extremals obtained in~\cite{el_max}:
\begin{align}
&J = \frac{\partial(x_t, y_t, \t_t)}{\partial(\p, k, r)} = 
\frac{1}{\sqrt r \cos (\t_t/2)} 
\frac{\partial(x_t, y_t, \sin (\t_t/2))}{\partial(\p, k, \sqrt r)} = 
- \frac{32}{(1-k^2)k^2 r^{3/2} \D^2} J_2, \label{JJ2} \\
&J_2 = c_2 z^2 + c_1 z + c_0, \qquad z = \tsp \in [0, 1], \label{J2=} \\
&p = \sqrt r t /(2 k), \qquad \tau = \sqrt r (2 \p + t/k)/2, 
\qquad \D = 1 - k^2 \ssp \tsp, \label{prtk} \\
&c_2 = k^4 \ss \cc x_1, \label{c2scx1} \\
&x_1 = 2 \cc \ss \E^3(p) + (\dd(3-6\ssp)-(2-k^2) p \cc \ss)\E^2(p) \nonumber \\
&\qquad + (\dd(k^2-2) p(1-2 \ssp) + \cc \ss (k^2 (2 p^2 - 1 + 6 \ssp) \nonumber \\
&\qquad - 2(2+p^2)))\E(p) + \dd(2k^2 \ccp \ssp + (1-k^2)p^2 (2\ssp-1)) \nonumber \\
&\qquad + p \cc \ss (2(2+p^2)-k^2(3+(3-k^2)p^2+(2-k^2)\ssp)), \nonumber  \\ 
&c_0 = - k f_2(p,k) \, x_2,  \label{c0fzx2}\\
&x_2 = \dd \E^2(p) - k^2 \cc \ss \E(p) - (1-k^2) p^2 \dd, \nonumber \\
&f_2(p,k) = 2 (\dd(2-k^2)p - 2 \E(p)) + k^2 \ss \cc)/k, \nonumber \\
&c_2 + c_1 + c_0 = (1-k^2) c_0. \label{c2c1c0c0}
\end{align} 

\subsection{Preliminary lemmas}
\begin{lemma}
\label{lem:c0<0N2}
For any $p > 0$ and $k \in (0, 1)$, we have $c_0 < 0$ and $c_0 + c_1 + c_2 < 0$.
\end{lemma}
\begin{proof}
In view of decomposition~\eq{c0fzx2} and  equality~\eq{c2c1c0c0}, it is enough to show that
\be{fz>0x2>0}
f_2(p,k) > 0, \quad x_2 > 0 \qquad \forall \ p > 0 \quad \forall \ k \in (0, 1).
\ee

We have 
$$
\left(\frac{f_2(p)}{\dd}\right)' = k^4 \frac{\ccp \ssp}{\ddp},
$$
this identity means that $f_2(p)/\dd$ increases w.r.t. variable $p$. But $f_2(0) = 0$, thus $f_2(p) > 0$ for all $p > 0$ and $k \in (0, 1)$.  

Further, from the equalities 
\begin{align*}
&\left(\frac{x_2(p)}{\dd \E(p)}\right)' = \frac{(1-k^2)(\E(p)-p \ddp)^2}{\ddp \E^2(p)},\\
&\E(p) - p \ddp = \frac 2 3 k^2 p^3 + o(p^3) \not\equiv 0,
\end{align*} 
it follows that $x_2(p)/(\dd \E(p))$ increases w.r.t. $p$. Then the asymptotics 
$$
x_2(p) = \frac{4}{45} (1-k^2)p^6 + o(p^6) > 0, \qquad p \to 0,
$$
 implies that $x_2 > 0$ for  all $p > 0$ and $k \in (0, 1)$. 

Inequalities~\eq{fz>0x2>0} are proved, and the statement of this lemma follows.
\end{proof}

\begin{lemma}
\label{lem:pconjN2pKn}
For any $n \in \N$, $k \in (0, 1)$, $z \in [0, 1]$, we have $J_2(Kn,z,k) < 0$.
\end{lemma}
\begin{proof}
Fix any $n$, $k$, $p = Kn$ according to the condition of this lemma. It follows from decomposition~\eq{c2scx1} that $c_2 = 0$. Thus the function $J_2(z)$ becomes linear:
$J_2(z) = c_1 z + c_0$, $z \in [0, 1]$. By virtue of Lemma~\ref{lem:c0<0N2}, this linear function is negative at the endpoints of the segment $z \in [0, 1]$:
$$
J_2(0) = c_0 < 0, 
\qquad 
J_2(1) = c_1 + c_0 = c_2 + c_1 + c_0 < 0,
$$
 thus it is negative on the whole segment $[0, 1]$ as well.
\end{proof}

\begin{lemma}
\label{lem:pconjN2k0}
For any $p_1 > 0$ there exists $\hk = \hk(p_1) > 0$ such that for all $k \in (0, \hk)$, $p \in (0, p_1)$, $z \in [0, 1]$ we have $J_2(p,z,k) < 0$.
\end{lemma}
\begin{proof}
In order to estimate the function $J_2$ for small $k$, we need the corresponding asymptotics as $k \to 0$:
\begin{align}
&c_0 = k^8 c_{00} + o(k^8), 
\qquad
c_1 = k^{10} c_{10} + o(k^{10}), 
\qquad
c_2 = k^{12} c_{20} + o(k^{12}), 
\label{c012k0} \\
&c_{00} = -c_{10} = -\frac{1}{1024} (4p - \sin 4p) \, c_{01}(p), \label{c00c01} \\
&c_{01} = 4p^2 -1 + \cos 4p + p \sin 4p, \nonumber \\
&c_{20} = \frac{1}{8192} \cos p \sin p \, c_{21}(p),  \nonumber \\
&c_{21} = - 3 \cos 2p - 48 p^2 \cos 2p + 3 \cos 6p + 42 p \sin 2p - 64 p^3 \sin 2p + 2p \sin 6p, \nonumber
\end{align}
and the asymptotics as $(p,k) \to (0,0)$:
\be{c02kp0}
c_0 = - \frac{4}{135} k^8 p^9 + o(k^8 p^9), 
\qquad
c_2 = \frac{4}{4725} k^{12} p^{11} + o(k^{12} p^{11}),
\ee
all these asymptotic expansions are obtained via Taylor expansions of Jacobi's functions, see~\cite{max3}.

(1) 
From the equalities 
$$
\left(\frac{c_{01}}{p}\right)' = \frac{2(\sin 2 p - 2 p \cos 2 p)^2}{p^2}, 
\qquad 
c_{01} = \frac{128}{45} p^6 + o(p^6),
$$ 
it follows that $c_{01}(p) > 0$ for $p > 0$, whence in view of decomposition~\eq{c00c01} we obtain that $c_{00}(p) < 0$ for all $p > 0$.

Fix an arbitrary number $p_1 > 0$.

(2)
Choose any $p_0 \in (0, p_1)$. We show that there exists $k_{01} = k_{01}(p_0, p_1) \in (0, 1)$ such that
\be{J2<0k01}
J_2(p,z,k) < 0 \qquad \forall \ p \in [p_0, p_1] \ \ \forall z \in [0, 1] \ \  \forall k \in (0, k_{01}).
\ee
Taking into account equalities~\eq{c012k0}, we obtain a Taylor expansion as $k \to 0$:
$$
J_2(p,z,k) = k^8 c_{00}(p) + \frac{k^{10}}{10!} \frac{\partial^{10} J_2}{\partial k^{10}} (p, z, \tk), \qquad p \in [p_0, p_1], \ z \in [0, 1], \ \tk \in (0, k).
$$
By continuity of the corresponding functions, we conclude that 
\begin{align*}
&c_{00}(p) < - m, \qquad m = m(p_0, p_1) > 0, \\
&\frac{1}{10!} \frac{\partial^{10} J_2}{\partial k^{10}} (p, z, \tk) < m_1, \qquad m_1 = m_1(p_0, p_1) > 0,
\end{align*} 
whence $J_2 < k^8 (-m + k^2 m_1) < 0$ for $k^2 < k^2_{01} = m/m_1 > 0$. Inequality~\eq{J2<0k01} follows.

(3)
From asymptotics~\eq{c02kp0} and equality~\eq{c2c1c0c0} we conclude that 
$$
J_2 = - \frac{4}{135} k^8 p^9 + o(k^8p^9), \qquad (p, k) \to 0.
$$ 
Thus
$$
\exists \ p_0' > 0 \ \exists k_0' > 0 \ \forall \ p \in (0, p_0'] \ \forall k \in (0, k_0') \ \forall \ z \in [0, 1] \qquad J_2(p,z,k) < 0.
$$

(4)
Take $p_0' \in (0, p_1)$ and $k_0' \in (0,1)$ according to item (3) of this proof. Find $k_{01} = k_{01}(p_0', p_1)$ according to item (2). Set $\hk (p_1) = \min (k_0', k_{01}) > 0$. Then for any $k \in (0, \hk(p_1))$ we get the following: if $p \in (0, p_0]$, then $J_2 < 0$ by item (3), and if $p \in [p_0, p_1]$, then $J_2 < 0$ by item (2). So the number $\hk(p_1)$ satisfies conditions of this lemma. 
\end{proof}

\subsection[Absence of conjugate points on non-inflectional elasticae]
{Absence of conjugate points \\ on non-inflectional elasticae}
\begin{theorem}
\label{th:N2conj}
If $\lam \in N_2$, then the normal extremal trajectory $q(t) = \Exp_t(\lam)$ does not contain conjugate points for $t > 0$.
\end{theorem}
\begin{proof}
In view of the symmetry $\map{i}{N_2^+}{N_2^-}$, see~\cite{el_max}, it is enough to consider the case $\lam \in N_2^+$. 

Denote $\lam^1 = \lam$. Fix any $n \in \N$ and prove that the trajectory 
$$
q^1(t) = \Exp_t(\lam^1), \qquad \lam^1 = (\f, k^1, r) \in N_2^+,
$$
does not contain conjugate points 
$t \in (0, t^1_1]$, $t^1_1 = 2 k^1 K(k^1) n/\sqrt r$.

Consider the family of extremal trajectories 
\begin{align*}
&\g^s = \{q^s(t) = \Exp_t(\lam^s)\mid t \in [0, t^s_1]\}, \\
&\lam^s = (\f, k^s, r) \in N_2^+, \quad 
t^s_1 = 2 k^s K(k^s) n/\sqrt r, \quad s \in [0, 1], 
\end{align*}
 where the covector $\lam^1 = (\f, k^1, r)$ is equal to $\lam$ given in formulation of this theorem, and covector $\lam^0 = (\f, k^0, r)$ will be chosen below in such a way that the parameter $k^0$ is sufficiently small. 

According to the Lemma~\ref{lem:pconjN2k0}, choose a number $\hk(p^1) \in (0,1)$ corresponding to the  number $p^1 = K(k^1) n$. Choose any $k^0 \in (0, \hk(p^1))$ and set $\lam^0 = (\f, k^0, r) \in N_2^+$.

By Lemma~\ref{lem:pconjN2k0}, for any $p \in (0, p^1]$ and any $z \in [0,1]$ we have $J_2(p,z,k^0) < 0$. By Lemma~\ref{lem:pconjN2pKn}, for any $z \in [0,1]$ and any $k \in [k_0, k_1]$ we have $J_2(K(k)n, z, k) < 0$.

Taking into account equality~\eq{JJ2} and relations~\eq{J2=}, \eq{prtk}, we conclude that the trajectory $\g^0$ does not have conjugate points at the segment $t \in (0, t^0_1]$, and for any trajectory $\g^s$, $s \in [0, 1]$, the endpoint $t = t^s_1$ is not conjugate. Now the statement of this theorem follows from Corollary~\ref{cor:homot_noconj}.
\end{proof}

\subsection{Absence of conjugate points for special cases}
The absence of conjugate points on for extremals $\lam_t \in N_2$ implies a similar fact for $\lam_t \in N_3 \cup N_6$.

\begin{theorem}
\label{th:N3N6conj}
If $\lam \in N_3 \cup N_6$, then the extremal trajectory $q(t) = \Exp_t(\lam)$ does not contain conjugate points for $t > 0$.
\end{theorem}
\begin{proof}
Let $\lam \in N_3 \cup N_6$. Since the set $N_3 \cup N_6$ belongs to the boundary of the domain $N_2$, one can construct a continuous curve $\map{\lam^s}{[0,1]}{N}$ such that $\lam^s \in N_2$ for $s \in [0, 1)$ and $\lam^1 = \lam$. 

Consider the family of extremal trajectories $q^s(t) = \Exp_t(\lam^s)$, $t > 0$, $s \in [0, 1]$. It follows from Th.~\ref{th:N2conj} that for $s \in [0, 1)$ the trajectory $q^s(t)$ does not contain conjugate points $t > 0$. Then we conclude from Corollary~\ref{cor:noconj_limit} that the trajectory $q^1(t) = \Exp_t(\lam)$ does not contain conjugate points for $t > 0$.
\end{proof}

\section{Final remarks}
\label{sec:final}
We sum up the study of conjugate points in Euler's elastic problem.

Any inflectional elastica contains an infinite sequence of isolated conjugate points. Visually, the first conjugate point occurs between the first and third inflection points. More intrinsically, the first conjugate point belongs to the interval $(T/2, 3T/2)$, where $T$ is the period of oscillation of the pendulum that parametrizes the vertical subsystem of the normal Hamiltonian system. The first conjugate point is contained in the arc bounded by the first Maxwell points corresponding to discrete symmetries~\cite{el_max}.

Non-inflectional and critical elasticae, circles and straight lines do not contain conjugate points. 

On the basis of this information about conjugate points and the description of Maxwell points obtained in~\cite{el_max}, one can study the global structure of the exponential mapping in Euler's elastic problem: describe the domains where the exponential mapping is diffeomorphic, and find a precise characterization of cut points. 
 Another interesting question for further study is the structure of the caustic in Euler's problem. 
This  will be the subject of our forthcoming work.

\bigskip\bigskip\bigskip\bigskip

\subsection*{Acknowledgment}
The author wishes to thank Professor A.A. Agrachev for proposing the problem and useful discussions during the work.

\newpage
\addcontentsline{toc}{section}{\listfigurename}
\listoffigures

\end{document}